\documentclass{amsart}
\usepackage{amssymb}
\ifx\pdfoutput\undefined
\usepackage{graphicx}  \DeclareGraphicsExtensions{.ps}

\else
\usepackage{amsfonts}
\usepackage[pdftex]{graphicx}  \DeclareGraphicsExtensions{.pdf,.mps}

\fi

\usepackage[T1]{fontenc}

\usepackage{epic, eepic}
\newcommand{\bc}{\mathbb C}

\newcommand{\im}{\operatorname{\rm Im}}
\newcommand{\mn}[1]{\Vert#1\Vert}
\newcommand{\op}{\operatorname{Op}}
\newcommand{\re}{\operatorname{\rm Re}}
\newcommand{\set}[1]{\left\{\,#1\,\right\}}
\newcommand{\w}[1]{\langle #1\rangle }

\def\squarebox#1{\hbox to #1{\hfill\vbox to #1{\vfill}}} 
\newcommand{\stopthm}{\hfill\hfill\vbox{\hrule\hbox{\vrule\squarebox 
                 {.667em}\vrule}\hrule}\smallskip} 

\newcommand{\F}{{\mathcal F}} 
\newcommand{\CC}{{\mathbb C}}
\newcommand{\CI}{{\mathcal C}^\infty }
\newcommand{\CIc}{{\mathcal C}^\infty_{\rm{c}} }
\newcommand{\CIb}{{\mathcal C}^\infty_{\rm{b}} }
\newcommand{\Op}{{\operatorname{Op}^{{w}}_h}}

\newcommand{\RR}{{\mathbb R}}

\newcommand{\ZZ}{{\mathbb Z}}
\newcommand{\NN}{{\mathbb N}}

\newcommand{\supp}{\operatorname{supp}}

\newcommand{\tr}{\operatorname{tr}}
\newcommand{\rest}{\!\!\restriction}

\renewcommand{\Re}{\mathop{\rm Re}\nolimits}
\renewcommand{\Im}{\mathop{\rm Im}\nolimits}

\theoremstyle{plain}

\newtheorem{thm}{Theorem}
\newtheorem{prop}{Proposition}[section]
\newtheorem{lem}[prop]{Lemma}

\theoremstyle{definition}

\numberwithin{equation}{section}

\title
{Pseudospectra of semi-classical (pseudo)differential operators}
\author[N. Dencker]{Nils Dencker}
\address{Matematikcentrum,  Lunds Universitet, Box 118, S-221 00 Lund, Sweden}
\email{dencker@maths.lth.se}
\author[J. Sj{\"o}strand]{Johannes Sj{\"o}strand}
\address{Centre de Math{\'e}matiques, {\'E}cole Polytechnique \\
UMR 7460, CNRS \\
F-91128 Palaiseau, France }
\email{johannes@math.polytechnique.fr}
\author[M. Zworski]{Maciej Zworski}
\address{Mathematics Department, University of California \\
Evans Hall, Berkeley, CA 94720, USA}
\email{zworski@math.berkeley.edu}

\begin{document}    
   
\maketitle   
   
\section{Introduction}   
\label{in}

The purpose 
of this note
is to show how some results from the
theory of partial differential equations apply to the study of
{\em pseudo-spectra} of non-self-adjoint operators which is a topic
of current interest in applied mathematics -- see \cite{Dav3} and \cite{Tre1}.

We will consider operators which arise from the quantization 
of bounded functions on the phase space $ T^* \RR^n $. 
For stronger results in the analytic case we will assume that our functions 
are holomorphic
and bounded in  tubular complex neighbourhoods of $ T^* \RR^n \subset 
\CC^{2n} $. 

Let us present the results in a typical example to which they apply.
We consider
\[ P( h )  = -h^2 \Delta + V ( x ) \,, \]
a semi-classical Schr\"odinger operator. 

We define the {\em semi-classical pseudospectrum} of the Schr\"odinger
operator $ P ( h ) $
as 
\[ \Lambda (p)  = \overline{ \{ \xi^2 + V ( x ) \; :\; 
( x, \xi) \in \RR^{2n} \,, \  \Im \langle \xi 
, V' ( x ) \rangle \neq 0 \} } \,, \]
noting that in the analytic case $ \Lambda (p) $ is either empty or
the closure of the set of all values of $ p= \xi^2 + V ( x ) $.

The following result (see \S \ref{sps}) shows that the resolvent is
large inside the pseudo-spectrum. We first state it
in the case of Schr\"odinger operators satisfying the assumptions above:

\begin{thm}
\label{thm00}
Suppose that $ P ( h ) = - h^2 \Delta + V ( x ) $, with
$ V \in \CI ( \RR^n ) $.

Then, there exists an open dense subset of $ \Lambda ( p ) $ such that
for any $ z $ in that subset there exists
$  u ( h ) \in L^2 ( \RR^n ) $ with the property 
\begin{equation}
\label{eq:t0}
  \| ( P ( h ) - z ) u ( h ) \| = {\mathcal O }( h^\infty  ) 
\| u ( h )\| \,. \end{equation}
In addition $ u ( h ) $ is localized to a point in phase
space, $ ( x , \xi ) $, with $ p ( x , \xi ) = z $.
More precisely,
$ WF_h ( u ) = \{ ( x , \xi ) \}$, where the wave front set,  $ WF_h ( u ) $,
is defined in \eqref{eq:wf}. Finally, for every compact $ K \Subset \Lambda ( p ) $
the above result holds uniformly for $ z \in K $ in the natural sense.
If the potential is real analytic 
then 
we can replace $ h^\infty $ by $ \exp ( - 1 / Ch ) $.
\end{thm}
This result was proved by Davies \cite{Dav2} for Schr\"odinger operators
in one dimension,
but as was pointed out in \cite{Zw}, it follows in great generality
from a simple
adaptation of the 
now classical results of H\"ormander \cite{Hor1} and Duistermaat-Sj\"ostrand 
\cite{DS}. The main point is that, unlike in the case of normal operators,
the resolvent can be large  on open sets as $ h \rightarrow 0 $.
That is particularly striking when $ P ( h ) $  has only discrete 
spectrum. 

To guarantee that we can for instance assume that
\begin{gather}
\label{eq:pot1}
\begin{gathered}
  | \partial^\alpha_x V ( x ) | \leq C_\alpha
 ( 1 + | x |)^{m-|\alpha| }  \,, \ 
\ \  
  ( 1 + |x |^m + |\xi|^2 ) 
/C  \leq | \xi^2 + V ( x ) | \,, \ \  |( x, \xi ) | \geq 
C \,.
\end{gathered}
\end{gather}
where $ m > 0$. This is the simplest example of the behaviour of the potential:
we can make weaker assumptions on $ V $ -- see the end of 
Sect.\ref{sps}. In the analytic case, we assume that \eqref{eq:pot1} holds as
 $ |x| \rightarrow \infty $, $ | \Im x | < c_0 $ (and we only need it with 
$ |\alpha| = 0 $).

The classical symbol $ p = \xi^2 + V ( x )  $
avoids all sufficiently negative values and the Fredholm theory
guarantees
that $ P ( h ) $
has discrete spectrum for $ h$ small enough (see \S \ref{rsq}).

We can, in place of the Schr\"odinger operator, $ P (h ) $, 
consider the operator with a bounded symbol, 
$ ( P ( h ) - z_1)^{-1} (P ( h ) - z_2)$, $ z_2 \neq z_1 $,
and this shows
that it is sufficient to consider quantization of bounded functions, with 
all derivatives bounded,
\[ p \in \CI_{\rm{b}} ( T^* \RR^n ) = \{ u \in \CI ( T^* \RR^n ) \; : \;
\forall \; \alpha \in \NN^n_0 \ \partial^\alpha u \in L^\infty ( T^* \RR^n ) 
\} \,.\]
In that case we give a more general definition of the semi-classical
pseudospectrum:
\begin{equation}
\label{eq:pdef}
\Lambda (p ) = \overline{ p ( \{ m \in T^* \RR^n \; : \; \{ p , \bar p\} (m) 
\neq 0 \}) } \,,
\end{equation}
where we used the Poisson bracket:
\[ \{ f , g \} = H_f g \,, \ \ 
H_f \stackrel{\rm{def}}{=} 
\sum_{ j=1}^n \partial_{\xi_j} f \partial_{ x_j}   - 
\partial_{x_j} f \partial_{ \xi_j}  \,.\]
The non-vanishing of $ \{ p , \bar p \} $ is a classical equivalent of 
the operator not being normal -- see \eqref{eq:prod1} and \eqref{eq:prod2} 
below. We note that in the analytic case we have
\[  \Lambda (p ) = \emptyset \ \text{ or } \ \Lambda ( p ) = \Sigma ( p ) 
\,,
\]
where we put
\[ \Sigma ( p ) = \overline{ p ( T^* \RR^n ) } \,.\]

In that more general setting we can restate our result as 
\begin{thm}
\label{thm0}
Suppose that $ n \geq 2 $, $ p \in \CI_{\rm{b}} ( T^* \RR^n ) $
and that $ p^{-1} ( z ) $ is compact for a dense set of
values $ z \in \CC$.
If $ P ( h ) $ has  the principal part given by 
$ p^w ( x , hD ) $ then the conclusions of Theorem \ref{thm00} hold.

If in addition $ p $  has a bounded holomorphic
continuation to $ \{ ( x , \xi ) \in \CC^{2n} \,, \ |\Im ( x , \xi ) | \leq 
1/ C \}$ then the conclusions of Theorem \ref{thm00} hold with 
$ h^\infty $ replaced by $ \exp ( - 1 / Ch ) $.

If $ n =1 $ then the same conclusion holds provided that the 
assumptions of Lemma \ref{l:2.2}$'$ are satisfied.

\end{thm}

We will see in \S \ref{dps} that, in general, we cannot construct an 
almost solution $ u ( h ) $ at an arbitrary interior
point of $ \Lambda ( p ) $, $ z $. 

In simple one dimensional examples we can already see that the spectrum,
$ \sigma ( P  ( h )) $, 
typically lies deep inside the pseudo-spectrum, $ \Lambda ( p ) $
 (the set of values of $ p$ in the analytic case) --  see \cite{Dav1},\cite{Dav2},\cite{Tre1}. 
Consider for instance the following 
non-self-adjoint operator $ P ( h ) = ( h D_x)^2 + i ( h D_x ) + x^2 $.
A formal conjugation 
$$ e^{ -x/ 2h } P ( h ) e^{ x / 2h }  = ( h D_x)^2 + x^2 + \frac14 \,, $$
shows that the spectrum of $ P ( h ) $ is given by $ ( 2n + 1) h + 1/4 $,
while 
\[ \Lambda ( p ) = \{ z \; : \;  \Re z \geq (\Im z)^2 \}\,, \ p = \xi^2 + 
i \xi + x ^2 \,.
\]

To see these phenomena for general operators we need to 
make assumptions on $ z_0 \in \partial \Lambda ( p ) $. The first one is 
the principal type condition,
\begin{equation}
\label{eq:prin}
p ( x, \xi ) = z_0 \ \Longrightarrow dp ( x, \xi 
 ) \neq 0 \,, \ \ m \in T^* \RR^n \,.
\end{equation}
Then we assume an exterior cone condition:

\begin{quote}
\begin{center}
There exists a truncated cone in $ \CC \setminus \Lambda ( p ) $ with vertex at $ z_0 $.
\end{center}
\end{quote}
More precisely,
\begin{equation}
\label{eq:cone}
\exists \; \epsilon_0 > 0 \,, \; \theta_0 \in \RR \ \ 
\text{ such that } \  ( z_0 + ( 0 , \epsilon_0 ) e^{ i ( \theta_0 -
\epsilon_0 , \theta_0 + \epsilon_0 ) } ) \cap \Lambda ( p ) = \emptyset \,,
\end{equation}
and a dynamical condition: if $ q = i e^{ - i \theta_0 } ( p - z_0 ) $, then
\begin{equation}
\label{eq:dyn}
\text{ No trajectory of $ H_{\Re q} $ can remain in $ q^{-1} ( 0 ) $ for 
an unbounded period of time.}
\end{equation}

Under these conditions we have the following 
\begin{thm}
\label{thm2}
Suppose that $ p \in \CI_{\rm{b}} ( T^* \RR^n ) $ and that  the principal part of  
$ P ( h )$ is given by $ p^w ( x , h D ) $.
If  $ z_0 \in \partial \Lambda ( p ) $ satisfies \eqref{eq:prin}, 
\eqref{eq:cone}, and \eqref{eq:dyn} then for any $ M > 0$, and 
for $ h < h_0 ( M )$, $ 0 < h_0  ( M ) $, 
\[   \{ z \; : \; | z - z_0 | < M h \log (1/h) \} \cap \sigma ( P ( h ) ) 
= \emptyset \,. \]

If in addition $ p $ 
 is a bounded holomorphic 
function in a complex tubular of $ \RR^n $.
then there exists $ C_0 > 0 $ 
such that
\[   \{ z \; : \; | z - z_0 | < 1/ C_0 \} \cap \sigma ( P ( h ) ) 
= \emptyset \,. \]
\end{thm}

In \S \ref{bn} we will show that if \eqref{eq:dyn} is violated then,
for a large class of 
{\em dissipative operators}, the spectrum lies arbitrarily close 
(as $ h \rightarrow 0 $) to the boundary of the pseudo-spectrum.

At the boundary of the pseudo-spectrum we may expect an improved bound
on the resolvent when some additional non-degeneracy is assumed. 
The result below is based on subelliptic estimates \cite[Chapter 27]{Horb}
and we borrow our notation from there.
If $p = p_1+ i p_2 \in C^\infty$ with real
valued $p_j$ then we define the repeated Poisson brackets
$$
p_I = H_{p_{i_1}} H_{p_{i_2}}\dots
H_{p_{i_{k-1}}}p_{i_k}
$$ 
where $I = (i_1,i_2,\dots,i_k) \in \set{1,2}^{k}$ 
and $|I| = k$ is the order of the bracket.

We say that $z_0 \in \partial{\Lambda}(p)$ is of {\em finite type} for ~$p$ if
\eqref{eq:prin} holds at $ z_0 $, 
$p^{-1}(z_0)$ is compact, and  
for any $(x_0,{\xi}_0) \in p^{-1}(z_0)$ there exists $k \ge 1$ and
$I \in \set{1,2}^{k}$ such that
\begin{equation}\label{finitetypecond}
 p_I(x_0,{\xi}_0) \ne 0.
\end{equation}
The {\em order} of $p$ at $w =(x_0,{\xi}_0)$ is 
\begin{equation}
 k(w) = \max \set{j \in \ZZ:\ p_I(w) = 0 \text{ for $|I| \le j$}}.
\end{equation}
The order of $z_0$ is the maximum of the
order of $p$ at $(x_0,{\xi}_0)$ for $(x_0,{\xi}_0) \in p^{-1}(z_0)$.
We say that $p$ satisfies condition $(P)$ if the imaginary part of
$qp$ does not change sign on the bicharacteristics of the real part of
$qp$, for any $0 \ne q \in C^\infty$.

As shown in \cite[Corollary 27.2.4]{Horb}, $ k ( w ) > k $ if and
only if
\begin{equation}
\label{eq:nonob}
 \forall \; z \in \CC\,, \ j \leq k \ \;  (H_{ \Re z p } )^j\Im z p 
( w ) = 0 \,, \end{equation}
and this provides a reformulation of the assumptions 
of the following
\begin{thm}\label{mainprop}
Assume that $p \in \CI_{\rm{b}} (T^* \RR^n ) $, and that the
principal part of ~$P(h)$ is ~$p^w ( x ,h D ) $. If $z_0 \in \partial {\Lambda}(p)$
is of finite type for ~$p$ of order $k \ge 1$, then $k$ is {\em even}
and for $ h < h_0 $, $ 0 < h_0 $, 
\begin{equation}
\label{eq:sube}
 \| ( P ( h ) - z_0 )^{-1} \| \leq C  h^{ - \frac{k}{k+1} } \,, \end{equation}
In particular, there exists $c_0 > 0$ such that
\begin{equation}
\set{z:\ |z-z_0|\le c_0h^{\frac{k}{k+1}}} \cap {\sigma}(P(h)) = \emptyset
\qquad\text{$0 < h \le h_0$}.
\end{equation}
\end{thm}
In one dimension this result was proved in \cite{Zwup}, and in 
some special cases by Boulton \cite{Bou} who also showed that the 
bounds are optimal. As was demonstrated by 
Trefethen \cite{Tr1} that is also easy to see numerically.

A simple higher dimensional example to which the theorem applies can 
be constructed as follows. Let $ W \in \CI_{\rm{b}} ( \RR^2 ) $ 
be a non-negative function, vanishing on the circle $ x_1^2 + x_2^2 = 1$.
Consider 
\[ P ( h ) = -h^2 \Delta + i W ( x ) + i ( x_1^2 + x_2^2 -1 )^m \,, \ \ 
\text{ with $ m$ even.} \]
Then the estimate \eqref{eq:sube} holds for $ z_0 > 0 $ uniformly on 
compact subsets of $ ( 0 , \infty ) $, with $ k = 2m $. The increase in 
$ k $ is due to the (simple) tangency of some bicharacteristics of the real
part to the set where the imaginary part vanishes. 

We conclude by pointing out that we could have defined the semi-classical 
pseudospectrum of $ P ( h ) $, $ {\Lambda} ( P ) $, as the closure of
the set of points $ z $ at which \eqref{eq:t0} holds. We have shown that
\[ \Sigma ( p ) \supset \Lambda ( P ) \supset \Lambda ( p ) \,.\]
An equality is not true in general but we could perhaps hope for 
\[ \overline{ \Lambda ( P ) ^\circ } = \Lambda ( p ) \,,\]
under suitable assumptions.

Another important topic not explored in this paper is the behaviour of
evolution operators $ \exp ( i t P / h ) $ for non-normal $ P $'s, and
its relation to semi-classical pseudospectra.

\medskip

\noindent

{\sc Acknowledgements.} The third author is
grateful to the 
National Science Foundation for partial support under the 
grant DMS-0200732. He would also like to thank Mike Christ and Nick
Trefethen for helpful discussion.

\section{Review of semi-classical quantization}
\label{rsq}

For simplicity of presentation we will consider the  case of
semi-classical quantization of functions
$ p \in \CI _{\rm{b}} ( T^* \RR^n ) $,
that is that $ p $ is bounded with bounded derivatives of all orders.

In the analytic case we will assume that 
 $ p ( x, \xi ) $ is 
bounded and holomorphic in a tubular neighbourhood of $ T^* \RR^n 
\simeq \RR^{2n } \subset \CC^{2n} $. 
As pointed out in the
introduction, the case of functions which omit a value in $ \CC $ and
which tend to infinity as $ ( x , \xi ) \rightarrow \infty $ can 
be reduced to this case (see also the remark at the end of \S \ref{sps}).

We use the Weyl quantization, 
\begin{equation}
\label{eq:2.weyl}
 p ^w ( x , h D_x ) u  = \frac{1}{ (2 \pi h )^n } 
\int \int  p \left ( \frac{ x+ y } 2 , \xi \right) e^{ \frac{i}{h} 
\langle x - y , \xi \rangle } u ( y ) d y d \xi \,.
\end{equation}
which for $ p \in \CI_{\rm{b}} ( T^* \RR^n ) $ gives operators bounded on 
$ L^2 ( \RR^n ) $ -- see \cite[Chapter 7]{DiSj}. We can consider more 
general operators,
\[ P ( h ) \sim \sum_{ j =0}^\infty h ^j p_j^w ( x , h D ) \,, \ \ \]
where in the case of analytic symbols we assume that 
\[ P(x,\xi ;h)\sim \sum_{j=0}^\infty h^j
p_j(x,\xi ) \,, \] 
in the space of bounded holomorphic functions in a tubular neighbourhood
of the real phase space.
Although it is not strictly speaking necessary for our final 
conclusions, in the analytic case we make an additional assumption that 
\begin{equation}
\label{eq:tube}
  | p_j ( z, \zeta ) | \leq C^j j^j \,, \ \ ( z , \zeta ) \in \CC^n \,, \ \ 
|\Im ( z , \zeta ) | \leq 1/C \,. \end{equation}
That allows us exponentially small errors in the expansions.

The product formula of the Weyl calculus says that
\begin{equation}
\label{eq:prod1}
p_1^w ( x , h D) \circ p_2 ^w ( x, h D) = 
( p_1 \sharp_h p_2 )^w ( x , hD; h) \,,\end{equation}
where 
\[
(p_1\# p_1)(x,\xi ;h)={e^{{ih\over 2}\omega ((D_x,D_\xi ),(D_y,D_\eta
))}p_1(x,\xi )p_2(y,\eta )_\vert }_{y=x,\,\eta =\xi } \]
has the following asymptotic expansion
\begin{equation} 
\label{eq:prod2}  p_1 \sharp_h p_2 ( x, \xi ; h ) \sim 
\sum_{ k=0}^\infty \frac{1}{k!} \left( \frac{ i h }{2} \omega ( 
( D_x, D_\xi ), (D_y , D_\eta)) \right)^k p_1 ( x , \xi ) p_2 ( y , 
\eta )  |_{ y = x, \eta = \xi } \,,
\end{equation} 
with $ \omega = \sum_{ j =1}^n d \xi_j \wedge d x_j $, the symplectic
form on $ T^* \RR^n $, and $ D_\bullet = (1/i) \partial_\bullet $. 
The expansion determines $ p_1 \sharp_h p_2 $ up to a term in 
 $ {\mathcal O} (h^\infty ) \CI_{\rm{b}} $. 
In the analytic case by summing up to $ k \sim 
1 / (C h) $ we can obtain $ {\mathcal O} (e^{ - 1/ ( Ch )  } ) $ errors -- 
see \cite{SjA}.

A basic tool of microlocal  analysis is the FBI transform:
\[  T : L^2 ( \RR^n ) \rightarrow L^2 ( T^* \RR^n ) \,, \]
defined by
\[ T u ( x , \xi ) = c_n h^{ - \frac{3n}{4}} \int_{\RR^n } 
e^{ \frac{i}{h} ( \langle x - y , \xi \rangle + i  | x - y |^2/2 ) }
u ( y ) d y \,. \]
Roughly speaking its r\^ole can be described as follows. The phase 
space properties of $ u \in L^2 ( \RR^n ) $ are reflected by the 
behaviour of $ T u \in L^2 ( T^* \RR^n ) $ as $ h \rightarrow 0 $. 
In this note we will only deal with $h$-dependent smooth functions 
with a tempered behaviour, $ | D^\alpha u | = {\mathcal O} ( h^{ - N_\alpha} )
$. The notion of the {\em wave front set} of $ u $, $ WF_h ( u ) $, 
explains the localization statement in Theorem \ref{thm0} (see also
Theorem \ref{thm0}$'$ below). 
In the $ \CI $ case the $h$-wavefront set is defined by 
\begin{equation}
\label{eq:wf}
  ( x^0 , \xi^0 ) \notin WF_h ( u )  \ \Longleftrightarrow  \ 
\forall \;
N \ 
| T u ( x , \xi ) | \leq C_N h^{-N}  \text{ for $ ( x , \xi) $ in a 
neighbourhood of $ ( x^0 , \xi^0 ) $,} \end{equation}
and in the analytic case by 
\begin{equation}
\label{eq:wf'}
  ( x^0 , \xi^0 ) \notin WF_h ( u )  \ \Longleftrightarrow \ 
\exists \;
c > 0 \ 
| T u ( x , \xi ) | \leq e^{ - c/ h } \text{ for $ ( x , \xi) $ in a 
neighbourhood of $ ( x^0 , \xi^0 ) $.} \end{equation}
In the $ \CI $ case we can also characterize $ WF_h ( u ) $ using 
pseudodifferential operators:
\[   ( x^0 , \xi^0 ) \notin WF_h ( u )  \ \Longleftrightarrow \
\exists \; p \in \CI_{\rm{c}} ( T^* \RR^n ) \,, \ p ( x^0, \xi^0 ) 
\neq 0 \,, \  p ( x , h  D) u = {\mathcal O} ( h^\infty ) \,.\]

In the analytic case we will need to understand the action of $ P ( h ) $
on microlocally weighted spaces $ H( \Lambda_{t G} ) $ whose definition,
in the simple setting needed here,
we will now recall -- see \cite{HeSj} for the origins of the method,
and \cite{mar} for a recent presentation. 

The complexification of the symplectic manifold $ T^* \RR^n $, $T^* \CC^n $
is equipped with the complex symplectic form, $ \omega_\CC $ and two 
natural real symplectic forms $ \Im \omega_\CC $ and $ \Re \omega_\CC $.
We see that $ T^* \RR^n $ is Lagrangian with respect to the first form
and symplectic with respect to the second one. In general we call a 
submanifold satisfying these two conditions an {\em IR-manifold}. 

Suppose that $ G \in \CI_{\rm{c}} ( T^* \RR^n ) $. We associate to it a
natural family of IR-manifolds:
\begin{equation}
\label{eq:2.ir}
\Lambda_{ t G } = \{ \rho + i t H_G ( \rho ) \; : \; \rho \in T^* \RR^n \}
\subset T^* \CC^n \,,  \ \ \text{with $ t \in \RR $ and $ |t | $ small.}
\end{equation}
Since $ \Im ( \zeta dz ) $ is closed on $ \Lambda_{t G } $, there exists
a function $ H_t $ on $ \Lambda_{ t G } $ such that
\[ d H_t = - \Im  ( \zeta dz )|_{ \Lambda_{t G } } \,,\]
and in fact we can write it down explicitely, parametrizing $ \Lambda_{tG} $ 
by $ T^* \RR^n $:
\[ H_t ( z, \zeta ) = - \langle \xi, t \nabla_\xi G ( x , \xi ) \rangle
+ t G ( x , \xi ) \,, \ \  ( z, \zeta ) = ( x , \xi ) + i t H_G ( x , \xi ) \,.  \]

The associated spaces $ H ( \Lambda_{ t G } ) $ are defined as follows. 
The FBI transform, $ T u ( x , \xi ) $, is analytic in $ ( x , \xi ) $ and we
can continue it to $ \Lambda_{tG} $. That defines $ T_{ \Lambda_{tG} }
u \in \CI ( \Lambda_{ t G} )  $.  Since $ \Lambda_{t G } $ differs from 
$ T^* \RR^n $ on a compact set only, $  T_{ \Lambda_{tG} } u $ is 
square integrable on $ \Lambda_{t G} $. 

The spaces $ H ( \Lambda_{t G} ) $ are defined by putting $ h$-dependent
norms on $ L^2 (\RR^n ) $:
\[ \| u \|_{ H( \Lambda_{ tG} ) } ^2 = \int_{ \Lambda_{t G} } 
| T_{\Lambda_{tG}} u ( z , \zeta )|^2 e^{ - 2 H_t ( z , \zeta ) / h } 
( \omega |_{\Lambda_{ t G} } )^n / n! \,. \]

The main result
relates the action of a pseudodifferential 
operator to the 
multiplication by its symbol. Suppose that $ p_1 $ and $ p_2 $ 
are bounded and holomorphic in 
a neighbourhood of $ T^* \RR^n $ in $ \CC^{2n} $ (see \eqref{eq:tube}). Then
for $ t $ small enough
\begin{equation}
\label{eq:cofe}
\begin{split} 
\langle p_1^w ( x , h D) u , p_2^w ( x, h D) v \rangle_{ H ( \Lambda_{ t G} ) } = &  
\langle (p_1|_{ \Lambda_{ t G } } ) T_{\Lambda_{tG}} u , 
 (p_2|_{ \Lambda_{ t G } } ) T_{\Lambda_{ tG}} v 
\rangle_{ L^2 ( \Lambda_{ t G } , e^{ - 2 H_t / h } (\omega |_{ \Lambda_{tG}
})^n/n! )} \\
& \ \ \ \ \ \ \ + \; {\mathcal O} ( h ) \| u \|_{ H ( \Lambda_{ t G} ) } 
\| v \|_{ H ( \Lambda_{ t G} ) } \,,
\end{split}
\end{equation}
see \cite{HeSj},\cite{mar}. In particular, 
by taking $ p_1 = p $ and $ p_2 = \bar p   $, and $ u = v $
we obtain
\begin{equation}
\label{eq:toep}
\| p^w ( x , h  D) u \|^2 _{ H ( \Lambda_{ t G} ) }  
= \|  p|_{ \Lambda_{ t G } }  T_{\Lambda_{tG}} u \|^2
_{ L^2 ( \Lambda_{ t G } , e^{ - 2 H_t / h } (\omega |_{ \Lambda_{tG}
})^n/n! )}  + {\mathcal O} ( h ) \| u \|^2 _{  H ( \Lambda_{ t G} ) }   
\,. 
\end{equation}

For the use in the next section 
we also recall some basic facts about {\em positive Lagrangian submanifolds}
of a complex symplectic manifold $ T^* \CC^n $. A complex plane
$ \lambda $, of (complex) dimension $ n$ is Lagrangian and positive if
\begin{equation}
\label{eq:2.pos}
\forall \; 
X, Y \in \lambda \ \ \
\omega_\CC ( X , Y ) = 0 \,, \ \ i \omega_\CC ( \bar X ,  X ) \geq 0 \,.
\end{equation}
The crucial characterization is given as follows
(see \cite[Proposition 21.5.9]{Horb}):
\begin{gather}
\label{eq:2.proj}
\begin{gathered}
\text{ $ \lambda \subset T^* \CC^n $ is a positive Lagrangian plane}
\ \Longleftrightarrow \ \lambda = \{ ( z , A z ) \; : \; z \in \CC^n \}
\\
\text{ where $ A = A_1 + i A_2 $ is a symmetric matrix 
with $ A_1 $ real, and $ A_2 $ positive definite.}
\end{gathered}
\end{gather}

\section{Semi-classical pseudo-spectrum}
\label{sps}

In \S \ref{in} we defined the semi-classical pseudospectrum, 
$ \Lambda ( p ) $,  as the closure of the set of values of $ p$.
We define some additional sets
\begin{equation}
\label{eq:ldef}
\begin{split}
& \Lambda_\pm ( p ) = \{ p ( x, \xi ) \; : \; 
\pm \{ \Re p , \Im p \} ( x, \xi ) > 0 \} \subset p ( T^* \RR^n ) \\
& \Sigma_{\infty} (p ) = \{ z \; : \; \exists \; ( x_j , \xi_j ) \rightarrow
\infty \  \  \lim_{ j \rightarrow \infty } p ( x_j , \xi_j ) = z \} \,,
\end{split}
\end{equation}
that is, $ \Sigma_\infty ( p ) $ is the set of limit
points of $ p $ at infinity.

In the $ \CI $ case Theorem \eqref{thm0} follows immediately from 
a semi-classical reformulation of the non-propagation of
singularities \cite{DS},\cite{Hor1},\cite{Hor2} -- see \cite[Section 26.3]{Horb}
and \cite{Zw}. 

The analytic case is also well known (see \cite{KK}) but since 
a ready-to-use reference is not available we include a proof. It can
also be adapted to give a self-contained proof in the $ \CI $ case.

\medskip
\noindent
{\bf Theorem \ref{thm0}$\bf '$.} {\it Suppose that $ n \geq 2 $ and
$ p ( x, \xi ) $ 
satisfies the assumptions in \S \ref{rsq} in the analytic 
case, and that $ \Lambda_- ( p ) $
is given by \eqref{eq:ldef}. 
Then
\[  \overline{\Lambda_- ( p )} \supset 
{ \Lambda ( p ) \setminus
\Sigma_\infty ( p ) } \,, \]
and for every 
$ z \in \Lambda_- ( p ) $, and every $(x^0,\xi ^0)\in T^*{\RR}^n$
with 
\[ p(x^0,\xi ^0)=z \,, \ \ \{ \Re p,\Im p\} (x^0,\xi ^0)<0\,, \]
 there exists
$0\ne u(h)\in L^2({\RR}^n)$ such that
such that
\begin{gather}
\label{eq:top'}
\begin{gathered}
 \| ( P ( h ) - z ) u ( h ) \| = {\mathcal O }( e^{ - 1/ Ch } ) 
\| u ( h )\| \,,  \ \ WF_h ( u ( h ) ) = \{ ( x^0 , \xi^0 ) \}  \,.
\end{gathered}
\end{gather}
If $ n =1 $ then the same conclusion holds provided that the 
assumptions of Lemma \ref{l:2.2}$'$ are satisfied.
}
In dimension one the theorem holds as well but further assumptions need
to be made on $ p$ -- see the remark after Lemma \ref{l:2.2}.

Before the proof we want to stress the need for an open dense subset
$ \Lambda_- ( p ) $.  
One could ask if any interior point of $ \Lambda ( p ) 
\setminus \Sigma_\infty ( p ) $
is an ``almost eigenvalue'' or ``quasimode'' in the sense of \eqref{eq:top'}.
That is not so as shown by

\renewcommand\thefootnote{\dag}

\medskip
\noindent
{\bf Example.}
Consider the following bounded analytic function on $ T^* \RR $:
\[ p ( x , \xi ) = \frac{ \xi^2 - 1 + i \xi x^2  ( 1 + x^2)^{-1}
 }{ 1 + \xi^2  + i \xi x^2 ( 1 + x^2)^{-1} } \,.\]
We see  that $ p^{-1} (0) = \{ ( 0 , 1 ) , ( 0 , -1 ) \}$, and that 
$ 0$ is an interior point of the pseudospectrum, 
$ 0 \in \Lambda( p )^\circ $. Also, $ 0 $ is a boundary point of images 
of neighbourhoods of $ ( 0 , \pm 1 ) $ under $ p$. 

 Near $ ( 0 , \pm 1 ) $, $ p $ is microlocally equivalent to a 
non-vanishing multiple of $ \xi  + i x^2 $. An explicit computation shows 
that 
the inverse of the models are bounded by $ h^{-\frac23} $, and a localization
argument\footnote{It is an easy one dimensional version of the
argument in \S \ref{dencker}.}  then shows that 
\[  \| p^w ( x , h D ) ^{-1} \|_{ L^2 \rightarrow L^2 }  \leq h^{-\frac23} 
\,. \]
Hence, $ 0 $ is not a quasi-mode. It should be stressed that the 
vanishing of the Poisson bracket $ \{ \Re p , \Im p \} $ (which occurs
in this example at $ p^{-1}(0) $) is not enough to guarantee the absence of
a quasi-mode. A violation of the condition $ \overline \Psi $  (see 
\cite[Sect.26.4]{Horb}) can produce quasi-modes with the simplest example
coming from adapting \cite[Theorem 26.3.6]{Horb} as in \cite{Zw}: 
$ p ( x , \xi ) = \xi - i x^k $, with $ k > 1$ odd.

\medskip

\renewcommand\thefootnote{\ddag}

We start the proof of Theorem \ref{thm0}$'$ with the discussion of
$ \Lambda_\pm ( p ) $. 
To establish that $ \Lambda_- ( p ) $ is dense 
we need the following result
of Melin-Sj\"ostrand \cite[Lemma 8.1]{MeSj}:
\begin{lem}
\label{l:2.1}
Suppose that $ n \geq 2 $ and
that $ d \Re p $,  $ d \Im p $ are linearly independent on 
$ p^{-1} ( z ) $. If $ \omega $ is the symplectic form on $ T^* \RR^n $
then 
\[ \{ \Re p , \Im p \} \lambda_{ p ,z} = \frac{\omega^{n-1}}{(n-1)!} \big|
_{ p^{-1} ( z ) } \,, \]
where $ \lambda_{ p , z } $ is the Liouville measure on $ p^{-1} ( z) $:
$ \lambda_{p ,z } \wedge d \Re p \wedge d \Im p = \omega^n / n! $. 

In particular, for any compact connected component of $ p ^{-1} ( z ) $, $ \Gamma$, 
we have 
\[ \int 
_\Gamma \{ \Re p , \Im p \} \lambda_{ p , z } ( d \rho ) = 0 \,. \]
\end{lem}

As an immediate consequence we see that $ \overline{\Lambda_- ( p )} = 
\Lambda ( p ) $ if the assumptions of Theorem \ref{thm0} are satisfied and that
in general we have the following
\begin{lem}
\label{l:2.2}
If the assumptions on $ p $ are satisfied, $ n \geq 2 $, and
either $ \Lambda_+ ( p ) $ or $ \Lambda_- ( p )  $ are non-empty, then 
$ \Lambda_+ ( p ) \cup \Lambda_- ( p ) $ is dense in $ \Lambda ( p )$, and
\[ \overline { \Lambda_{ \pm } ( p ) } \supset  \Lambda ( p )  \setminus 
\Sigma_\infty 
( p ) \,. \]
\end{lem}
\begin{proof}
Assume that $\{ \Re p,\Im p\}\not\equiv 0$. Then $H \stackrel{\rm{def}}{=}
\{ \rho \in T^*{\RR}^n
;\, \{\Re p,\Im p\} (\rho )=0\}$ is an analytic hypersurface without any
interior
points. Consequently, every value $z=p(\rho )$ with $\rho \in H$ can be
approximated by values $z_j=p(\rho _j)$ with ${\RR}^{2n}\setminus H\ni \rho
_j\to \rho $, and $\Lambda _+(p)\cup \Lambda _-(p)$ is open and dense in
$p(T^*{\RR}^n)$.

Since under our assumption 
 $ (p ( T^* \RR^n ))^\circ \neq \emptyset $, an elementary version of 
the Morse-Sard theorem implies that $ d \Re p $ and $ d \Im p $ are independent
on $ p^{-1} ( z) $ for $ z $ in a dense open set 
$ \Omega \subset  \Lambda ( p ) \setminus \Sigma_\infty ( p )  $. 
Lemma \ref{l:2.1} then shows
that $ \Lambda_+ ( p ) \cap \Omega = \Lambda_- ( p ) \cap \Omega $, completing
the proof of  the lemma.
\end{proof}

\medskip
\noindent
{\bf Remark.} 
In the case of dimension one a different argument, based on 
elementary topological considerations, is needed and some assumptions
have to be made on $ p $. To see that consider for instance
\[ p ( x , \xi ) = \frac{ ( \xi + i x )^2 }{ 1 + x^2 + \xi^2 } \,, \ \ 
 \{ \Re p , \Im p \} ( x, \xi ) > 0 \,, \ ( x , \xi ) \neq ( 0 , 0 ) \,. \]
For $ p $'s arising from Schr\"odinger operators
considered in \S \ref{in}
we always have
\begin{equation}
\label{eq:sgn}  \sum_{ m \in p^{-1} ( z ) } {\rm{sgn}}\;  \{ \Re p , \Im p \} ( m ) = 0 \,,
\end{equation}
for a dense set of values $ z $. In fact,
$ p ( x , \xi ) = p ( x , - \xi ) = z $, and the set of values $ z $
corresponding to $ \xi \neq 0 $ is dense in the set of values for which
the bracket is non-zero. Now we simply notice that
\[  \{ \Re p , \Im p \} ( ( x , \xi )  ) = - 
 \{ \Re p , \Im p \} ( ( x , - \xi )  )  \,. \]

\medskip

\noindent
{\bf Lemma \ref{l:2.2}$\bf '$.} {\it Suppose that $ n = 1 $ and
in addition to the general assumptions, each component of 
$ \CC \setminus \Sigma_\infty ( p ) $ has a non-empty 
intersection with $ \complement \Lambda ( p ) $.
Then the conclusions of Lemma \ref{l:2.2} and \eqref{eq:sgn} (for 
a dense set of values) hold.}
\begin{proof}
Let $\Omega$ be a
component of $ \CC \setminus \Sigma_\infty (p)$. Then 
\[ \iota \stackrel{\rm{def}}{=} {\rm var\,\,
arg}_{\gamma(z)} (p-z) \]
is independent of $z\in \Omega$ if $\gamma (z)$ is the
positively oriented circle $|(x,\xi )|=R(z)$ with $R(z)$ large enough. 
For $z\in
\Omega\setminus\Lambda(p)$ $\iota$ is zero (for mapping degree reasons) 
and hence
it is zero for all $z\in \Omega$. If $z\in \Lambda (p)\cap\Omega$ is a regular
value, we get 
\[0=\iota = 2\pi  \sum_{m\in p^{-1}(z)}  {\rm{sgn}}\; \{ \Re p,\Im p\} 
(m) \,, \]
so $z$ belongs to both $\Lambda_+ ( p ) $ and $\Lambda_- ( p ) $.
\end{proof}

In the remainder of the proof there is no restriction on the dimension.

\medskip
\noindent
{\em Proof of Theorem \ref{thm0}$'$:}
We can assume that $ z = 0 $ and 
we follow the now standard procedure of the {\em complex WKB construction}
associated to a positive Lagrangian submanifold of the complexification of 
$ T^ * \RR^n $. We start with the geometric construction of that 
submanifold. Since $ \{ \Re p , \Im p \} \neq 0 $ we have $ d_\xi p \neq 0 $,
and we can assume that $ \partial_{\xi_1} p (  x^0 , \xi^0 ) \neq 0 $. Let 
$ \phi_0 $, be a real analytic function defined in a neighbourhood of 
$ y^0 \in \RR^{n-1} $, $ x^0 = ( x_1^0 , y^0 ) $, with the properties
\[  \phi_0 ( y^0) = 0 \,, \ \ d \phi_0 ( y^0 ) = \eta^0 \,, \ \xi^0 = 
( \xi_1^0, \eta^0 ) \,, \ \ \Im d^2 \phi_0 ( y^0 ) \gg 0 \,, \]
where the Hessian $ d^2 \Im \phi_0 $ is well defined at $ y^0 $ as 
$ d \Im \phi_0 = 0 $. We will make further assumptions on $ \phi_0 $
later.

We then define $ \Lambda_0 \subset T^* \CC^n $, locally near $ 
( x^0 , \xi^0 )$, as follows
\[ \Lambda_0 = \{ ( x_1^0 , y ; \xi_1 ( y ) , d_y \phi_0 ( y ) ) \; : \; 
p ( x_1^0 , y ; \xi_1 ( y ) , d_y \phi_0 (y ) ) = 0\,, \  \xi_1^0 ( y^0 ) 
= \xi_1^0 \} \,, \]
where we know that the function $ \xi_1 ( y ) $ is locally defined and
analytic from our condition $ \partial_{\xi_1 } p \neq 0 $. 
Using holomorphic continuation we obtain a locally defined submanifold of
$ T^* \CC^n $,   $ \Lambda_0 \cap T^* \RR^n = \{ ( x^0, \xi^0 ) \}$. 
This submanifold is isotropic with respect to the complex
symplectic form and its tangent spaces\footnote{These tangent spaces,
$ T_\rho \Lambda_0 $, are complex linear subspaces of $ T_\rho \CC^n $.}
are positive, in the sense that \eqref{eq:2.pos} is satisfied without
the condition on the dimension. 

\renewcommand\thefootnote{$\star$}

For $ t \in \CC $, $ |t | < \epsilon $, the complex flow, $ \Phi_t $,
exists by the Cauchy-Kovalevskaya Theorem:
\begin{gather*}
\Phi_t ( z  , \zeta  ) = ( z ( t ) , \zeta ( t ) ) \,, \ \ 
z ( 0 ) = z \,, \ \zeta ( 0 ) = \zeta \,, \\
z' ( t) = \partial_\zeta p ( z ( t ) , \zeta ( t ) ) \,, \ \
\zeta' ( t ) = - \partial_z p ( z ( t) , \zeta ( t) ) \,,
\end{gather*}
that is $ d/dt ( z ( t ) , \zeta ( t ) ) = H_p ( z ( t) , \zeta ( t ) ) $, 
$ \omega_\CC  ( \bullet , H_p ) = d p $. 

We then define
\[ \Lambda = \bigcup_{ t \in \CC \,, |t | < \epsilon }
\Phi_t ( \Lambda_0 ) \subset T^* \CC^n \,, \]
which is Lagrangian with respect to $ \omega_\CC $. 

We now want to guarantee that
tangent spaces to $ \Lambda $ are positive in the sense of
\eqref{eq:2.pos}.
We first note that
\begin{gather*}
i \omega_\CC (  \overline{ t H_p } , t  H_p )
= i|t|^2  \{ \bar p, p \} = - 2 \{ \Re p , \Im p \} |t|^2 >  \gamma |t|^2
> 0  \,,
\end{gather*}
by the assumptions of the theorem. We
have
\begin{gather*}
T_{ ( x^0 , \xi^0 )} \Lambda = T_{ ( x^0 , \xi^0 )} \Lambda_0
+ {\rm{span}}_\CC H_p ( x^0 , \xi^0 ) \,, \\
 i \omega_\CC (  \overline{  X  + t H_p }  , X  + t H_p )
=  i \omega_\CC (  \overline X  , X ) +
2 \Im ( t d p^* ( X ) )
+ |t|^2 i \omega_\CC (  \overline{  H_p } ,  H_p ) \,,
\end{gather*}
where $ p^* ( \rho ) = \overline {p ( \bar \rho ) } $. Hence, for
positivity, we need to
show that we can choose $ \phi_0 $ so that for $ X \in
T_{ ( x^0 , \xi^0 )} \Lambda_0 $,
\[  | d p^* ( X ) |^2 <  \alpha i  \omega_\CC (  \overline X  , X )
\,, \ \  \alpha = -2 \{ \Re p , \Im p \} ( x^0 , \xi^0 ) \,. \]
A calculation in local coordinates $ ( z_1 , z' ; \zeta_1, \zeta' ) $
shows that this follows from
\begin{gather*} \| A + \Phi B \|^2 <  \alpha \min {\rm{Spec}} ( \Im \Phi )
\,,
\\
\Phi = \phi_0'' \,, \ \  A = |p_{\zeta_1}'|^{-1}
i \left( p_{\zeta_1} \overline{ p}_{z'} - \overline{p}_{\zeta_1}
{ p_{z'}} \right)
\,, \ \
B = |p_{\zeta_1}'|^{-1}
i \left( p'_{\zeta_1} \overline{ p}_{\zeta'} -
\overline{p}_{\zeta_1}{ p_{\zeta'}}
\right) \,. \end{gather*}
The vectors $ A $ and $ B $ are real, and hence we can choose the complex
matrix $ \Phi $ so that $ A + ( \Re \Phi ) B = 0 $. This leaves us with
\[ \| ( \Im \Phi ) B \| ^2 <
 \alpha \min {\rm{Spec}} ( \Im \Phi )  \,, \]
which can be arranged by making $ \Im \Phi $ sufficiently small\footnote{An
alternative, and slicker, way of proceeding is by first observing that
the positivity is invariant under affine linear canonical 
transformations. Using
that, and a multiplication by a non-vanishing factor, we can assume that
$ p = \xi_n - i x_n + {\mathcal O} ( ( x, \xi)^2 ) $ and that $ (x^0 , \xi^0 ) 
= ( 0, 0 ) $. It is then straightforward to find $ \phi ( x ) $ with the
desired properties.}.

From \eqref{eq:2.proj} we see that $ \Lambda  \subset \{ p = 0 \} $
is locally a graph, and
since it is Lagrangian, a graph of a differential of a phase function 
$ \phi $. Since the tangent plane is positive
\eqref{eq:2.proj} shows that  the Hessian of that phase function
has a positive definite imaginary part:
\[ \Lambda = \{ ( z , d_z \phi ( z) ) \} \,, \ \ p ( z , d_z \phi ) = 0 \,, 
\ \phi ( x^0 ) = 0\,, \ d_z \phi ( x^0 ) = \xi^0 \,, \ \Im d^2_z \phi ( 
x^0 ) \gg 0 \,. \]
We also note that $ \Lambda \cap T^* \RR^n = \{ ( x^0 , \xi^0 ) \} $ which 
corresponds to the fact that $ \Im d_z \phi \neq 0 $ for $ z \neq x^0 $.

Once the phase function has been constructed we apply the usual 
WKB construction:
\[ v ( z , h ) \sim e^{ i \phi ( z ) /h } \sum_{ j =0 } ^\infty 
a_j ( z) h^j \,,  \]
where we will want 
the coefficients, $ a_j $'s to be  holomorphic near $ z = x^0 $, and to satisfy bounds
$ | a_j (  z) | \leq C^j j^j $. They are constructed so that 
\[\left(  \sum_{ j < 1/ Ch } 
p_j^w ( z , h D_z ) h^j \right)
\left(  e^{ i \phi ( z ) /h } \sum_{ j < 1/ Ch  }  
a_j ( z) h^j \right) = {\mathcal O} ( e^{ - 1/ Ch } ) \,.\]
Here $ p^w $ denotes the Weyl quantization of a holomorphic symbol $ p ( z , 
\zeta ) $ acting on holomorphic functions (compare to \eqref{eq:2.weyl}):
\[ p ^w ( z , h D_z ) u  = \frac{1}{ (2 \pi h )^n } 
{\int \! \! \! \int}
_{\Gamma_z} p \left ( \frac{ z+ w } 2 , \zeta \right) e^{ \frac{i}{h} 
\langle z - w , \zeta \rangle } u ( w ) d w d \zeta \,,
\]
where the contour $ \Gamma_z $ is suitably chosen  -- see \cite[Section 4]{SjA} for a discussion of the  general case.

The transport equations for $ a_j $'s then are:
\[ \sum_{ k = 1}^n \partial_{\zeta_k} 
p_0 ( z, d_z \phi ( z) ) \partial_{z_k }
a_j ( z) + i p_1 ( z, d_z \phi ( z)  ) a_j  = A_j ( z ) \,, \]
where $ A_j ( z) $ depends on $ a_l $'s with $ l < j   $, and we put $ 
a_0  (x^0_1 , y)  = 1 $.
It is now classical that the solutions 
satisfy $ |a _j | \leq C^j j^j $ near $ x^0 $ -- see \cite[Theorem 9.3]{SjA}. 
The real quasi-mode is
obtained by restricting to the real axis and by truncating $ v (z , h ) $:
\[ u ( x, h ) = \chi ( x ) v ( x , h ) \,, \ \  \chi ( x)  =  1\,, \ \ 
| x - x^0 | < \delta \,, \ \ \supp \chi \subset {\mathbb B} ( x^0 , 2 \delta) 
\]
where $ \delta $ is small. Since the construction has shown that 
$ \Im \phi \geq | x - x^0 |^2 / C $, the cut-off function $ \chi $ 
does not destroy the exponential smallness of the  error. 
\stopthm

For completeness, and later use in \S \ref{bn}, we include a result on the
discreteness of the spectrum.

\begin{prop}
\label{p:3}
Suppose that $ p \in \CI_{\rm{b}} ( T^* \RR^n ) $. Let $ \Omega $ be 
an open connected ($h$-independent) set, satisfying
\[  \overline \Omega \cap \Sigma_\infty ( p ) = \emptyset\,, \ \ 
 \Omega \cap  \complement \Sigma ( p ) \neq \emptyset \,.\]
Then 
$ ( p^w ( x , h D) - z )^{-1} $, $ 0 < h < h_0 ( \Omega) $, $ z \in \Omega $,
is a meromorphic family of operators 
with poles of finite rank.

 In particular, for $ h$ sufficiently small, the spectrum of 
$ p^w ( x , h D) $ is discrete in any such set.
\end{prop}
\begin{proof}
If  $ \Omega $ satisfies the assumptions of the propostion then 
there exists $ C > 0$ such that for every $ z \in \Omega $, we have
$ | p ( x , \xi ) -  z| > 1/ C $ if $ | ( x , \xi ) | > C $. 
The assumption that $ \Omega \cap \complement \Sigma ( p ) \neq \empty$ implies
that for some $ z_0 \in \Omega $, $ ( p ( x, \xi ) - z_0 )^{-1} \in \CI_{\rm{b}} ( T^* \RR^n ) 
$. Let $ \chi \in \CIc ( T^* \RR^n ; [0,1]) $ be equal to $ 1 $ in a sufficiently large 
bounded domain. The remarks above show that
\[  r ( x , \xi ; z ) = \chi ( x , \xi ) ( z_0 - p ( x , \xi ) )^{-1} 
+ ( 1 - \chi ( x , \xi )) ( z - p ( x , \xi )^{-1} \,,\]
is in $ \CIb ( T^* \RR^n ) $. The symbol calculus reviewed in \S \ref{rsq} 
then gives
\[  \begin{split}
& r^w ( x , h D, z ) ( z - p^w ( x, h D )  ) = I + {\mathcal O}_{L^2 \rightarrow L^2 }
( h ) + K_1 ( z ) \,, \\ 
& ( z - p^w ( x, h D )  ) r^w ( x , h D, z )  = I + {\mathcal O}_{L^2 \rightarrow L^2 }
( h ) + K_2 ( z )\,,\end{split} \]
where $ K_j ( z ) $, $ j = 1, 2 $ are compact operators on $ L^2 ( \RR^n ) $,
depending holomorphically on $z$ and vanishing for $z=z_0$.

By the analytic Fredholm theory we conclude that $ ( z - p^w ( x, h D ) )^{-1} 
$  is 
meromorphic in $ \Omega $ for $ h $ sufficiently small.
\end{proof}

\medskip
\noindent
{\bf Remark.} The same result holds for $ P ( h) $ of the form considered
in \S \ref{rsq} with $ p_j \in \CI_{\rm{b}} ( T^* \RR^n ) $: the lower 
order terms do not affect the meromorphy when $ h $ is small. We also 
comment on the case presented in \S \ref{in}.

Suppose that $ m ( x, \xi  ) $ is an admissible weight function, that is
a positive function on $ T^* \RR^n \simeq \RR^{2n} $ satisfying
\[ \forall \; X, Y \in \RR^{2n}  \ \ 1 \leq  m ( X )  \leq C \langle X - Y \rangle^N
m ( Y ) \,, \]
for some fixed $ C  $ and $ N $. Following \cite{DiSj} for symbols 
satisfying $ |\partial_X^\alpha p ( X ) | \leq C_\alpha m ( X ) $ we can
define operator $ P = p^w ( x , h D ) $. In the analytic case we 
require that
\[  | p ( X ) | \leq m ( \Re X ) \,, \ \ |\Im X | \leq 1/C \,.\]
In the example given in \S \ref{in} we can take $ m ( x , \xi ) 
= \langle \xi\rangle^2 + \langle x \rangle^m  $. 

Under an ellipticity assumption 
\[  | p ( X ) | \geq m ( \Re  X ) /C  \,, \ 
\ |X | \geq C \,, \ \ |\Im X | \leq 1/C \,,\]
we obtain an invertibility: if $ z_1 \not \in \overline{p ( \RR^{2n})} $
then  $ P - z_1 $ is invertible. If we define the operator
\[  Q = ( P - z_1)^{-1} ( P - z_3 )\,, \  \ z_3 \neq z_1 \,,\]
 then the  
resolvents of $ Q $ and $ P $ are related by 
\[  ( Q - \zeta)^{-1} =  ( 1 - \zeta)^{-1} ( P - z_1 ) \left( P - 
\frac{ \zeta z_1 - z_3 }{ \zeta - 1 } \right)^{-1} \,, \]
so that the reduction of Schr\"odinger operators to the case of operators
with bounded symbols was justified. 


\section{Decrease of pseudo-spectrum by a change of norms}

\label{dps}

In this section we will prove Theorem \ref{thm2}. That will be done
by using a dynamically defined function $ G$ which grows on 
the bicharacteristics of $ \Re p$.

In the analytic case 
we will now use the microlocally weighted spaces, $ H ( \Lambda_{t G } ) $,
the construction of which was recalled in \S \ref{rsq}, to 
decrease the pseudo-spectrum by changing the norm on $ L^2 ( \RR^n )$ 
in an $h$-dependent way. In particular this will show that under
assumptions \eqref{eq:prin}, \eqref{eq:cone}, and \eqref{eq:dyn}, 
the spectrum is separated from the boundary of the pseudo-spectrum.

\subsection{Construction of microlocal weights}

We start with 
\begin{lem}
\label{p:1}
If \eqref{eq:prin} and \eqref{eq:cone} hold,  then for $ q = 
 i e^{ - i \theta_0 } ( p - z_0 ) $
and  $\rho _0\in T^* \RR^n $ with $q (\rho _0)=0$, we have 
\begin{equation}
\label{eq:3}
\Im q\le 0  \ \text{ on } \ (\Re q)^{-1}(0) \cap \text{\em neigh}(\rho_
0, \RR^{2n}).
\end{equation}
Furthermore, $d\Re q(\rho _0)\ne 0$ and $d\Im q(\rho _0)=\lambda_0
d\Re q(\rho _0)$ for some $\lambda_0 \in {\RR}$.
\end{lem}

\medskip
\noindent
{\bf Notation.}
From now on we will assume, without loss of generality, that $ p = q $,
and $ z_0 = 0 $. We recall the 
conditions \eqref{eq:cone} and \eqref{eq:dyn}:
\begin{gather}
\label{eq:cone1} 
\exists \; \epsilon_0 > 0 
\text{ such that } \   ( 0 , \epsilon_0 ) i e^{ i ( -
\epsilon_0 ,  \epsilon_0 ) }  \cap \Lambda ( p ) = \emptyset \,,\\
\label{eq:dyn1} 
\text{ No trajectory of $ H_{\Re p} $ can remain in $ p^{-1} ( 0 ) $ for 
an unbounded period of time.}
\end{gather}

\medskip

\begin{proof}
We first recall that $d p(\rho _0)\ne 0$ by \eqref{eq:prin}.
If this differential were purely imaginary, it is
easy to see that we get a contradiction with \eqref{eq:cone1}. 
Hence $d\Re p(\rho _0)\ne 0$. Using \eqref{eq:cone1} again
it is clear that $ q \leq 0 $ on $ \Re q ^{-1} ( 0 ) $, near 
$ \rho_0 $, and that $ d \Im p ( \rho_0 ) $ cannot 
be independent of  $ d \Re p ( \rho_0 ) $.
\end{proof}

The next lemma gives a construction of the weight:
\begin{lem}
\label{l:4.2}
Suppose that \eqref{eq:prin},
\eqref{eq:dyn1} hold.
Then there exists $G\in \CI_{\rm{c}} (\RR^{2n};\RR)$, such that 
$H_{\Re p}G(\rho )>0$ for every $\rho \in
p^{-1}(0)$.
\end{lem}
\begin{proof}
The assumption \eqref{eq:dyn1}, gives a seemingly stronger statement:
\begin{equation}
\label{eq:dyn'}
\exists \; T_0 > 0 \ \forall \; (x, \xi)  \in ( \Re q )^{-1} ( 0 ) \ \exists \;0<  t < T_0 
\ \ \text{ such that } \ ( \Im q ) ( \exp ( t H_{ \Re q } ) ( ( x, \xi ) 
 ) ) \neq 0 \,,
\end{equation}
and it allows us a construction of a global weight $ G$.
Indeed, we first construct $G$ locally: Let $\gamma (t)$, $0\le t\le
T_1$ ($0\le T_1<T_0$) be a maximal $H_{\Re p}$ integral curve in
$p^{-1}(0)$. Then we can find a real-valued $G\in C_0^\infty $ with
support in a small neighbourhood of the image of $\gamma $ such that $H_{\Re
p}G\ge 0$ on $p^{-1}(0)$ with strict inequality on the image of $\gamma $.
We then get the $G$ of the lemma, by taking a finite sum of such local
$G$'s.
\end{proof}

\subsection{The $\CI $ case}

\def\cint{{1\over 2\pi i}\int}
\def\iint{\int\hskip -2mm\int}
\def\iiint{\int\hskip -2mm\int\hskip -2mm\int}
\def\buildover#1#2{\buildrel#1\over#2}
\font \mittel=cmbx10 scaled \magstep1
\font \gross=cmbx10 scaled \magstep2
\font \klo=cmsl8
\font\liten=cmr10 at 8pt
\font\stor=cmr10 at 12pt
\font\Stor=cmbx10 at 14pt
\def\aby{arbitrary}
\def\ably{arbitrarily}
\def\asy{asymptotic}
\def\bdd{bounded}
\def\bdy{boundary}
\def\coef{coefficient}
\def\coeff{coefficient}
\def\const{constant}
\def\Const{Constant}
\def\canform{canonical transformation}
\def\coef{coefficient}
\def\coeff{coefficient}
\def\cont{continous}
\def\diff{diffeomorphism}
\def\diffeo{diffeomorphism}
\def\de{differential equation}
\def\dop{differential operator}
\def\ef{eigenfunction}
\def\ev{eigenvalue}
\def\e{equation}
\def\eq{equation}
\def\fy{family}
\def\fu{function}
\def\F{Fourier}
\def\fop{Fourier integral operator}
\def\fourior{Fourier integral operator}
\def\fouriors{Fourier integral operators }
\def\hol{holomorphic}
\def\hm{homogeneous}
\def\indep{independent}
\def\lhs{left hand side}
\def\mfld{manifold}
\def\ml{microlocal}
\def\neigh{neighborhood}
\def\nondeg{non-degenerate}
\def\op{operator}
\def\og{orthogonal}
\def\pb{problem}
\def\Pb{Problem}
\def\pde{partial differential equation}
\def\pe{periodic}
\def\per{periodic}
\def\pert{perturbation}
\def\Prop{Proposition}
\def\pol{polynomial}
\def\pop{pseudodifferential operator}
\def\pseudor{pseudodifferential operator}
\def\res{resonance}
\def\rhs{right hand side}
\def\sa{selfadjoint}
\def\sc{semiclassical}
\def\schr{Schr\"odinger operator}
\def\sop{Schr\"odinger operator}
\def\st{strictly}
\def\stpsh{\st{} plurisubharmonic}
\def\strans{^\sigma \hskip -2pt}
\def\suf{sufficient}
\def\sufly{sufficiently}
\def\tf{transformation}
\def\Th{Theorem}
\def\th{theorem}
\def\tf{transform}
\def\trans{^t\hskip -2pt}
\def\top{Toeplitz operator}
\def\uf{uniform}
\def\ufly{uniformly}
\def\vf{vector field}
\def\wrt{with respect to}
\def\Op{{\rm Op\,}}

 Let
\begin{equation}
\label{4}
{C_1h\le \epsilon \le C_2h\log {1\over h},}
\end{equation}
where $C_1>0$ is large enough. We shall derive an estimate for the \e{}
\begin{equation*}
{P_\epsilon ( h ) u=v,}
\end{equation*}
when ${\rm WF}_h(u)$ is contained in a small \neigh{} of $p^{-1}(0)$, and
where 
\begin{equation*}
P_\epsilon ( h )  \stackrel{\rm{def}}{=} 
e^{\epsilon G/h}P ( h ) e^{-\epsilon G/h}=e^{{\epsilon \over h}{\rm
ad}_G}P( h ) \sim \sum_0^\infty  {\epsilon ^k\over k!}({1\over h}{\rm
ad}_G)^k(P ( h ) )\,, \  G = G^w ( x , h D) \,.  \end{equation*}
We note that the assumption on $ \epsilon $ and 
the boundedness of  $ {\rm ad}_G/h $ show that the expansion makes sense.
The operators $ \exp ( \epsilon G/ h ) $ are pseudo-differential in 
an exotic class $ S_\delta^{C_2}   $ for any $ \delta > 0 $ (see 
\cite{DiSj}) but that is not relevant here. We are essentially using 
the method of pseudodifferential operators of variable order \cite{Unt}.
In a related context of the absence of resonances similar methods
where recently used in \cite{Mar1}.

Dropping $ h$ in $ P ( h ) $ and 
using the same letters for \op{}s and and the corresponding symbols, we see
that
\begin{equation*}
{P_\epsilon =P+i\epsilon \{ p,G\} +{\mathcal O}(\epsilon ^2)=p+i\epsilon \{
p,G\} +{\mathcal O}(h+\epsilon ^2),}
\end{equation*}
so that
\begin{equation*}
\begin{split}
& {\Re P_\epsilon =\Re p-\epsilon \{ \Im p,G\} +{\mathcal O}(h+\epsilon ^2),}\\
& {\Im P_\epsilon =\Im p+\epsilon \{ \Re p,G\} +{\mathcal O}(h+\epsilon ^2).}
\end{split}
\end{equation*}
The positivity assumption  implies that
\begin{equation}
\label{9}
{\Im p+\lambda (x,\xi )\Re p\ge 0, }
\end{equation}
for a suitable smooth and real-valued function $\lambda $. Let us consider
\begin{equation}
\label{10}
{
\Im P_\epsilon +\lambda \Re P_\epsilon =
}
{\Im p+\lambda (x,\xi )\Re p+\epsilon \{ \Re p-\lambda \Im p,G\} +{\mathcal
O}(h+\epsilon ^2),}
\end{equation}
and let us 
write \eqref{9} as $\Im p=-\lambda \Re p+r$, where $r\ge 0$ vanishes on
$p^{-1}(0)$, so that $\nabla r$ also vanishes there. Then
\begin{equation*}
\begin{split}
\{ \Re p-\lambda \Im p,G\}  = & \{ (1+\lambda ^2)\Re p,G\} -\{ \lambda r,G\}= 
\\
& {(1+\lambda ^2)\{ \Re p,G\} +\{\lambda ^2,G\} \Re p -\{ \lambda r,G\} \ge
1/C,}
\end{split}
\end{equation*}
 since $\Re p$, $\{ \lambda r,G\}$ are small near $p^{-1}(0)$.

Using this and \eqref{9} in \eqref{10}, we get
\begin{equation}
\label{12}
{
\Im P_\epsilon +\lambda \Re P_\epsilon \ge {\epsilon \over C}+{\mathcal
O}(h+\epsilon ^2). }
\end{equation}

Consider
\begin{equation*}
{J=({1\over 2i}(P_\epsilon -P_\epsilon ^*)u\vert u)+({1\over 2}(\lambda
P_\epsilon +P_\epsilon ^*\lambda )u\vert u).}
\end{equation*}
On the one hand,
\begin{equation*}
{J=\Im (P_\epsilon u\vert u)+\Re (\lambda P_\epsilon u\vert u),}
\end{equation*}
so
\begin{equation}
\label{15}
{
\vert J\vert \le C\Vert P_\epsilon u\Vert \, \Vert u\Vert .
}
\end{equation}
On the other hand, the symbol of
$${1\over 2i}(P_\epsilon -P_\epsilon ^*)+{1\over 2}(\lambda P_\epsilon
+P_\epsilon ^*\lambda )$$
is equal to ${\mathcal O}(h)$ plus the expression \eqref{12}, so \eqref{12}
and the sharp G\aa{}rding inequality (see \cite[Theorem 7.12]{DiSj} 
or apply \eqref{eq:cofe} with $ G=0 $ and $ p_1 =  \Im P_\epsilon +
\lambda \Re P_\epsilon $, $ p_2 = 1 $) imply
\begin{equation*}
{J\ge {\epsilon \over 2C}\Vert u\Vert ^2.}
\end{equation*}
Combining this with \eqref{15}, we get
$${
\Vert u\Vert \le {2C^2\over \epsilon }\Vert P_\epsilon u\Vert
}$$
and in view of \eqref{4}, where $C_2>0$ can be \ably{} large, we conclude that
for every $C>0$, we have
$$
{
D(0,Ch\log {1\over h})\cap\sigma (P)=\emptyset ,\ 0<h<h(C),
}$$
when $h(C)>0$ is \sufly{} small.

\subsection{The analytic case}

To apply the theory of weighted spaces reviewed in \S \ref{rsq} 
we need one more

\begin{lem}
\label{l:4.3}
Suppose that \eqref{eq:prin}, \eqref{eq:cone1}, and \eqref{eq:dyn1} 
hold, and that $ G $ is given by Lemma \ref{l:4.2}. In the notation
of \eqref{eq:2.ir} we have, for sufficiently small $ t >0 $,
\begin{equation}
\label{eq:l.4} 
| p \rest_{ \Lambda_{ t G} } | >  t / C \,.
\end{equation}
\end{lem}
\begin{proof}
The proof of Lemma \ref{p:1} shows that if $\rho _0$ is a point with
$p(\rho _0)=0$ and we define $\lambda _0$, by 
$d\Im p(\rho _0)=\lambda _0 d\Re p(\rho _0)$, 
then for every $\epsilon >0$, there is a neighbourhood $W$
of $\rho _0$, such that if $\rho \in W$, then
\begin{equation}
\label{eq:5}
{
\Im p(\rho )\le \lambda _0 \Re p(\rho )+\epsilon | \Re p(\rho )| .
}
\end{equation}
Now,
$$p(\rho +itH_G)=p(\rho )-itH_pG(\rho )+{\mathcal O}(t^2)=p(\rho )-itH_{\Re
p}G+tH_{\Im p}G+{\mathcal O}(t^2).$$
It follows that
\[ \begin{split}
&(\Im p-\lambda _0 \Re p)(\rho +itH_G(\rho ))\\
&= (\Im p-\lambda
_0 \Re p)(\rho )-tH_{\Re p}G(\rho )-\lambda _0tH_{\Im p}G(\rho )+{\mathcal
O}(t^2) \\
&= (\Im p-\lambda _0\Re p)(\rho )-(1+\lambda _0^2)tH_{\Re p}G(\rho_0 )+
{\mathcal
O}(t\vert \rho -\rho _0\vert )+{\mathcal O}(t^2).
\end{split}
\]
Since $H_{\Re p}G(\rho _0)>0$,  for $\vert \rho -\rho _0\vert $ small
enough,
\begin{equation}
\label{eq:6}
{(\Im p -\lambda _0 \Re p-\epsilon \vert \Re p\vert )(\rho +itH_G(\rho
))\le} { -{\frac1C}t+{\mathcal O}(t\vert \rho -\rho _0\vert )+{\mathcal
O}(\epsilon t)\le -{\frac{1} {2C}}t,}
\end{equation}
where we also used \eqref{eq:5}. 
\end{proof}

As a consequence of the last lemma we see 
that $0$ is a removable point of
the pseudospectrum. To see this
we recall that  $H(\Lambda _{tG})$ is equal to $L^2$ as a space and that
the norms are equivalent for every fixed $h$ (but not uniformly with respect 
to $h$). 
The spectrum of $ P ( h ) $ therefore does not depend on whether 
we realize this operator on $L^2$ or on $H(\Lambda _{tG})$. 
We conclude that $0$ has an $h$-independent
 neighbourhood which is disjoint from the spectrum of $ P ( h ) $, when
$h$ is small enough. All this is summarized in the following result which 
is an immediate consequence of Lemma \ref{l:4.3} and \eqref{eq:toep}:

\medskip

\noindent
{\bf Theorem \ref{thm2}$'$.} 
{\em 
Suppose that $ z_0 \in \partial \Lambda ( p ) \setminus \Sigma_\infty ( p )  $ 
and that \eqref{eq:prin}, \eqref{eq:cone}, and \eqref{eq:dyn} hold.
In the notation of Lemma \ref{l:4.2} let us introduce introduce
the IR-manifold $\Lambda _{tG}=\{ \rho +itH_G(\rho );\, \rho 
\in{\RR}^{2n}\}$ for $t>0$ small enough. If
\begin{gather*}
 P ( h ) \sim \sum_{ j } h^j p_j ^2 ( x , h D ) \,, \  \ p_0 = p \,, \ 
\text{ $p_j$'s satisfy the assumptions of \S \ref{rsq},}
\end{gather*}
then 
\begin{gather*}
P( h ) - z_0 \; : \; H(\Lambda _{tG}) \ \longrightarrow \  H(\Lambda _{tG})\,, 
\end{gather*}  
has a bounded inverse for $ h $ small enough.
 In particular, for $ \delta $ small enough but independent of $ h $, 
\[  \sigma ( P ( h ) ) \cap D ( z_0 , \delta ) = \emptyset \,, \ \ 0 < h < h_0
\,. \]
}

\section{Proof of Theorem \ref{mainprop}}
\label{dencker}

We first observe that
that property of being of finite type and the order of the symbol is
invariant under multiplication with elliptic factors.  In fact, if $0
\ne q \in C^\infty$ then the repeated Poisson brackets of $\re qp$ and
$\im qp$ of order $\le j$ is a linear combination with smooth
coefficients of those of $p_1$ and $p_2$ and vice versa (see
Section~27.2 in \cite{Horb}).

Since we rely heavily on the localization argument based on Weyl calculus
of pseudodifferential operators we will follow \cite[Section~18.5]{Horb}
and introduce 
\begin{equation}\label{ghdef}
 g_h(dx,d{\xi}) = |dx|^2 + h^2|d{\xi}|^2 \,.
\end{equation}
We then 
find that $p(x,h{\xi}) \in S(1,g_h)$. We will also introduce 
the notion of {\em local pseudo-spectrum} defined as follows:
\begin{equation*}
 {\Lambda}(p({\Omega})) = \overline{p\left({\Omega} \cap \{ \{ p , 
\bar p \} \neq 0 \} \}  \right)}.
\end{equation*}

\subsection{Reduction to normal form}\label{reduct}

Assume that $z_0 \in \partial {\Lambda}(p)$ is a principal value
that is, \eqref{eq:prin} holds.
By subtracting $z_0$ we may assume that $z_0 = 0$.  Let $w_0 =
(x_0,{\xi}_0) \in p^{-1}(0)$, then since $|dp(w_0)| \ne 0$ we may
assume that $\partial_x p =0$ and $\partial_{\xi_j}p =0$, $j> 1$, at
~$w_0$ by a symplectic change of coordinates. By using the Malgrange
Preparation Theorem, we obtain
\begin{equation}\label{reducform}
p = q({\xi}_1 + i f(x,{\xi}')) = q\widetilde p\qquad {\xi}= ({\xi}_1,{\xi}')
\end{equation}
in a neighborhood ${\Omega}$ of $w_0$. Here $0 \ne q \in C^\infty$ and
$f(x,{\xi}') \in C^\infty$ is real valued. In the following, we put
$B_r = \set{w \in T^*\RR^n: \ |w-w_0| \le r}$ for $r > 0$.

\begin{lem}\label{boundarylem}
We have that $0 \in \partial{\Lambda}\left(p(B_r) \right)$
for sufficiently small $r >0$
if and only if $\pm f \ge 0$ on $B_r$ for a choice of sign, which
implies that $\widetilde p$ and ~$p$ satisfy condition ~$(P)$  on ~$B_r$. 
\end{lem}

Thus the property of a principal value being a boundary value is
preserved under multiplication with elliptic factors.

\begin{proof}
We find that $\pm f \ge 0$ in $B_r$ if and only if $\pm \arg
\left[\widetilde p(B_r)\setminus 0 \right] \subseteq [0, {\pi}]$
(otherwise it is ${} = [-{\pi},{\pi}]$).  Since $q(w)/q(w_0) \to 1$
for $w \in B_r$ when $r \to 0$, we find that $|\arg(q(w)/q(w_0))| <
{\pi}/2$ when $w \in B_r$ for small enough $r$. We find that $\pm
f \ge 0$ in ~$B_r$ if and only if $p(B_r)\setminus 0 = q\widetilde
p(B_r)\setminus 0 $ is contained in a sector $\set{z\in \bc:\ \arg z
\notin [a,b]}$ where $a \ne b$, and thus $0 \in \partial
{\Lambda}\left(p({B_{r}})\right)$ for small enough $r$.  In fact, if $f
\gtrless 0$ in $B_r$ for some $r>0$ then it is easy to find a
piecewise linear closed curve ${\gamma}_r$ in ~$B_r$ such that
$\widetilde p({\gamma}_r)$ (and thus $p({\gamma}_r)$) has winding
number equal to one in $\bc$. If ${\Pi}$ is a piecewise linear surface
containing ${\gamma}_\varrho$ then $p$ maps the bounded component of
${\Pi}\setminus {\gamma}_\varrho$ onto a neighborhood of the origin,
since $p(w_0) = 0$.
\end{proof}

We note that the lemma implies that \eqref{eq:cone} holds locally.

\begin{proof}[Proof of Theorem~\ref{mainprop}] 
As before, we may obtain $p$ of the form ~\eqref{reducform} near $w_0
= (x_0,{\xi}_0) \in p^{-1}(0)$ after subtracting a constant.  Since $0
\in \partial{\Lambda}(p)$ we find from Lemma~\ref{boundarylem} that $f
\ge 0$ after possibly switching $x_1$ to $-x_1$ and multiplying with
$-1$.  By the invariance we find that $\widetilde p_I(w_0) \ne 0$ for
some $I$ such that $|I| =k+1$ where $k$ is the order of $p$ at ~$w_0$,
which is less or equal to the order of $z_0 =0$.  Since $f \ge 0$ we
find that $H_f = 0$ when $f = 0$. Thus, the only non-vanishing
repeated Poisson bracket of order $k+1$ at ~$w_0$ is equal to
$\partial_{x_1}^{k}f(0,{\xi}_0')$, ${\xi}_0 = (0,{\xi}_0')$, and
since $f \ge 0$ we find that $k= 2l$ is even and
\begin{equation}\label{adcond}
\partial^{2l}_{x_1} f(0,{\xi}_0') > 0 
\end{equation} 
see \cite[Propositions 27.2.1, 27.2.2, and (27.3.1)]{Horb}.  
It can be seen from the (highly non-obvious) equivalent characterization
\eqref{eq:nonob}.

By choosing a suitable
cut-off function ${\psi}(x',{\xi}') \in C_0^\infty(T^*\RR^{n-1})$
supported near $(0,{\xi}_0')$ such that $0 \le {\psi}(x',{\xi}') \le
1$ and replacing $f(x,{\xi}')$ by ${\psi}(x',{\xi}')f(x,{\xi}')
\linebreak[0] + x_1^{k}(1- {\psi}(x',{\xi}'))$ we may assume that
$f\in C^\infty$ is uniformly bounded and
\begin{equation}\label{subcond}
 \partial^{k}_{x_1} f(x,{\xi}') >0\qquad \text{when $|x_1| \le c$.}
\end{equation}
By a change of variables we may assume that $c = 1$.

Since $p^{-1}(0)$ is compact, we may take $0 \le {\phi}_j \in
C_0^\infty(T^*\RR^n)$, $j = 1, \dots, N$, and ${\phi}_0 = 1 - \sum_{1
\le j \le N}^{}{\phi}_j$, such that $p$ is on the
form~\eqref{reducform} with $f$ satisfying~\eqref{subcond} in $\supp
{\phi}_j$, $j > 0$, and $|p| \ge c_0 > 0 $ on $\supp {\phi}_0$. Then we
obtain that
\begin{equation*}
 \mn {{\phi}_0^w(x,hD_x)u} \le C(\mn{p^w(x,hD_x)u} + h\mn{u})
\end{equation*}
by using that $p^{-1}(x,h{\xi}) \in S(1,g_h)$ on $\supp {\phi}_0$.
We find from Proposition~\ref{modelest} and~\eqref{subcond} that 
\begin{equation*}
 \mn {{\phi}_j^w(x,hD_x)u} \le Ch^{-\frac{k}{k+1}}(\mn{p^w(x,hD_x)u} +
 h\mn{u})\qquad j >0
\end{equation*}
since we may assume that $c = 1$ in ~\eqref{subcond}, and
\begin{equation}
\widetilde p^w(x,hD_x) {\phi}_j^w(x,hD_x) \cong
{\phi}_j^w(x,hD_x)(q^{-1})^w(x,hD_x)p^w(x,hD_x)
  \qquad\text{}
\end{equation}
modulo $\Op S(h, g_h)$. This gives 
\begin{equation*}
 \mn {u} \le \sum_{j =0}^N \mn{{\phi}^w_j(x,hD_x)u} \le
 C_0h^{-\frac{k}{k+1}}(\mn{p^w(x,hD_x)u} + h\mn{u}).
\end{equation*}
For small enough $h > 0$ we find that 
\begin{equation}
\mn {u} \le  C_0h^{-\frac{k}{k+1}}\mn{p^w(x,hD_x)u} 
\end{equation}
which proves Theorem~\ref{mainprop}.
\end{proof}

The estimates for the localized operators will be proved in the next subsection.

\subsection{The Model Operator}

In this section we shall consider the following subelliptic model operator
\begin{equation}\label{modelop} 
P(h) = hD_t + if^w(t,x,hD_x)\qquad 0 \le f \in C^\infty \qquad 0 < h
\le 1 
\end{equation}
which we assume satisfies ~\eqref{subcond}.
In this subsection we rely heavily on \cite[Section 27.3]{Horb}. 

\begin{prop}\label{modelest}
Assume $P(h)$ is given by ~\eqref{modelop}, where $f \in C^\infty$ is
uniformly bounded and  satisfies   
\begin{equation}\label{modeladcond}
\max_{j \le k}|\partial^j_{t} f(t,x,{\xi})| \ne 0  
\end{equation}
for $|t| \le 1$. Then we obtain 
\begin{equation}\label{propest}
\mn{u} \le C h^{-\frac{k}{k+1}} \mn{P(h)u}
\end{equation}
for $u \in C_0^\infty$ having support when $|t| \le 1$. 
\end{prop}

\begin{proof} 
By the non-negativity
of $f \in C^\infty$ we obtain from  \cite[Lemma~7.7.2]{Horb1} that 
\begin{equation}\label{fderest}
|\partial_x f|^2 +  |\partial_{\xi}f|^2 \le C f\qquad |t| \le 1.
\end{equation}
In the following, we shall assume that $|t| \le 1$.

Next, we write $P(h) = h P_h$ where 
\begin{equation}
P_h = D_t + i h^{-1}f^w(t,x,hD_x) = D_t + iF_h^w(t,x,D_x)
\end{equation}
with $F_h(t,x,{\xi}) =  h^{-1}f(t,x,h{\xi})\in C^\infty(\RR,
S(h^{-1}, g_h))$.
We obtain from ~\eqref{modeladcond} that
\begin{equation}
\max_{0 \le j\le k}|\partial_t^j F_h|^{\frac{1}{j+1}} \ge ch^{-\frac{1}{k+1}} 
\end{equation}
since $h^{-\frac{1}{j+1}} \ge h^{-\frac{1}{k+1}} $ when $0 \le j \le k$ and
$0 < h \le 1$.

Now we let $F_{h,w}(t)= F_h(t,w)$ for $w = (x,{\xi})$. 
By Lemma~27.3.4 in ~\cite{Horb} we obtain
\begin{equation}\label{subest}
h^{-\frac{1}{k+1}} \mn{u} + \mn{F_{h,w}(t)u}\le C\mn{(D_t + i
  F_{h,w}(t))u} \qquad \forall\,w
\end{equation}
for small enough $h>0$ and $ u \in C_0^\infty$ having support where
$|t| \le 1$. In fact,  
$$M_1(t) = \max_{j \le k} |\partial_t^j F_{h,w}(t)| \ge \max(|
F_{h,w}(t)|, ch^{-\frac{1}{k+1}}) \gg 1$$
when $0 < h \ll 1$ and $|t| \le 1$.

Next we shall introduce a new symbol class adapted to $F_h$. By
\eqref{fderest} we obtain that 
\begin{equation}
h |\partial_xF_h|^2 + h^{-1} |\partial_{\xi}F_h|^2 \le C|F_h|,
\end{equation}
which means that
\begin{equation}\label{dfest}
|F_h|^{g_h}_1 \le C_1  \sqrt{F_h}h^{-1/2}.
\end{equation}
Let 
\begin{equation}
m(t,w) = F_h(t,w) + h^{-1/k+1} \ge 1
\end{equation}
and $g_{m(t,w)} = g_h/hm(t,w)$, with $g_h$ given by ~\eqref{ghdef}. Then it
follows from ~\eqref{dfest} that $g_m$ is ${\sigma}$ ~temperate,
$g_m/g_m^{\sigma} = m^{-2} \le 1$ and that $m$ is a weight for
$g_m$. In fact,
$$ |F_h(t,w)-F_h(t,w_0)| \le C{\varrho}m(t,w_0) \qquad\text{when
$g_{m(t,w_0)}(w-w_0) \le {\varrho}^2,$}$$ 
since then $g_h(w-w_0)
\le {\varrho}^2hm(t,w_0)$. This gives the slow variation, and  since
$$g^{\sigma}_{m(t,w_1)}(w_1-w_0) = m(t,w_1)g_h(w_1-w_0)/h \ge
cm(t,w_0)$$ 
when $g_{m(t,w_0)}(w_1-w_0)
\ge c$, we obtain that $g_m$ is ${\sigma}$ ~temperate. 
It also follows from ~\eqref{dfest} that
$F_h \in S(m,g_m)$, since $|F_h|^{g_m}_1 =
\sqrt{m}h^{1/2}|F_h|^{g_h}_1$ and $$|F_h|^{g_m}_j \le C_j
(mh)^{j/2}|F_h|^{g_h}_j \le C_j m\,, \ j \ge 2 \,. $$

Next, we shall localize the estimate. Take a partition of unity
$\set{{\phi}_j(w)} \in S(1,g_{\varepsilon})$, where $g_{\varepsilon} =
h^{2{\varepsilon}-1}g_h$ and $0 < {\varepsilon} < \frac{1}{2k+2}<
\frac{1}{2}$ is fixed. This can be done uniformly in $h$, and
in the following all estimates are uniform in $h$.
Observe that since $2{\varepsilon} < \frac{1}{k+1}$ we
find $g_m \le g_{\varepsilon}$ for any $t$.  We assume that
${\phi}_j$ is supported in a sufficiently small $g_{\varepsilon}$
neighborhood of $w_j$, so that $m(t,w) \cong m(t,w_j)$ only varies
with a fixed factor in $\supp {\phi}_j$ for any $t$.  Since $\sum
{\phi}_j^2 = 1$ and $g_{\varepsilon} =
h^{4{\varepsilon}}g_{\varepsilon}^{\sigma}$ the calculus gives
\begin{equation}
\sum_j \mn{{\phi}_j^w(x,D_x)u}^2 -C h^{4{\varepsilon}}\mn{u}^2 \le
\mn{u}^2 \le \sum_j \mn{{\phi}_j^w(x,D_x)u}^2 + C
h^{4{\varepsilon}}\mn{u}^2
\end{equation}
for $u \in C_0^\infty$, thus for small enough $h$ we find 
\begin{equation}\label{pouest}
\sum_j \mn{{\phi}_j^w(x,D_x)u}^2  \le 2\mn{u}^2 \le
4 \sum_j \mn{{\phi}_j^w(x,D_x)u}^2\qquad\text{for $u \in C_0^\infty$}.
\end{equation}
We find from ~\eqref{subest} for small enough
~$h$ that
\begin{equation}\label{locest}
h^{-\frac{1}{k+1}} \mn{{\phi}_j^w(x,D_x)u} +
\mn{F_{h,w_j}(t){\phi}_j^w(x,D_x) u}\le C\mn{{\phi}_j^w(x,D_x)(D_t + i
  F_{h,w_j}(t))u} \qquad \forall\, j
\end{equation}
and the calculus  gives
\begin{equation}\label{errest}
{\phi}_j^w(x,D_x)\left(F_h^w(t,x,D_x) - F_{0,w_j}(t)\right)=R_j^w
(t,x,D_x) \in
S(m^{1/2}h^{-{\varepsilon}}, g_{\varepsilon}) 
\end{equation}  
since 
$$|F_h(t,w) - F_{h,w_j}(t)| \le
C\sqrt{m}h^{-{\varepsilon}}  \text{in $\supp {\phi}_j$} $$
 by ~\eqref{dfest}. 
In fact, we find that $|F_h|^{g_{\varepsilon}}_1 = h^{\frac{1}{2}
  -{\varepsilon}} |F_h|^{g_h}_1 \le Cm^{1/2}h^{-{\varepsilon}}$, and 
for $k \ge 2$ we have 
$$|F_h|^{g_{\varepsilon}}_k = h^{k(\frac{1}{2} 
  -{\varepsilon})} |F_h|^{g_h}_k \le C_kh^{\frac{k-2}{2} -
  k{\varepsilon}} \le C_k h^{-2{\varepsilon}} \le 
C_k\sqrt{m}h^{-{\varepsilon}}$$
 since ${\varepsilon}< \frac{1}{2k+2} < \frac{1}{2}$. 

We obtain from~\eqref{pouest} that
\begin{equation}
 \mn{F_h^w(t,x,D_x)u}^2 \le
4\sum_j \mn{F_{h,w_j}(t){\phi}_j^w(x,D_x)u}^2 + 4\w{R^w(t,x,D_x)u,u}  \qquad u
  \in C_0^\infty
\end{equation}
where $R^w(t,x,D_x) = \sum_j \overline R_j^w(t,x,D_x) R_j^w(t,x,D_x) \in \Op
S(m^{}h^{-2{\varepsilon}}, g_{\varepsilon})$. Since  
$$m = F_h +
h^{-1/k+1} > 0 \,, $$ 
 we obtain that $m^{-1} \in S(m^{-1},
g_m)$.  In fact, the  Leibniz' rule gives
\begin{equation*}
 \left|{m}^{-1}\right|^{g_m}_k \le C_k\sum_{r_1 + \dots + r_j
   =k} m^{-1-j}\prod_{i=1 }^j
 |m|^{g_m}_{r_i} \le C'_k m^{-1}.
\end{equation*}
Thus we obtain that
\begin{equation}
\left(m^{-1}\right)^w(t,x,D_x) m^w(t,x,D_x) \cong 1
\end{equation}
modulo $\Op S(m^{-2}, g_m) \subseteq \Op
S(h^{4{\varepsilon}},g_{\varepsilon})$, since $m^{-2} \le h^{4{\varepsilon}}$.
This gives
\begin{equation*}
 \mn{R^w u} \le C(h^{-2{\varepsilon}}\mn{(F_h +
   h^{-1/k+1})^w u} + h^{4{\varepsilon}}\mn {R^w u}) 
\end{equation*}
since $h^{2{\varepsilon}}R^w\left(m^{-1}\right)^w \in \Op S(1,
g_{\varepsilon})$.
We obtain by Cauchy-Schwarz that 
\[ \begin{split}
 |\w{R^wu,u} | & \le C h^{-2{\varepsilon}}\mn{(F_h^w +
   h^{-1/k+1})u}\mn{u} \\
& \le {\varrho}(\mn{F_h^wu}^2 +
 h^{-\frac{2}{k+1}}\mn{u}^2)  +
 C_{\varrho}h^{-4{\varepsilon}}\mn{u}^2
\end{split} 
\] 
for $u \in C_0^\infty$. We 
also find that there exists $H_{\varrho}> 0$ so that 
\begin{equation}
C_{\varrho} h^{-4{\varepsilon}}\mn{u}^2 \le {\varrho}
h^{-\frac{2}{k+1}}\mn{u}^2\qquad\text{when $0 < h \le H_{\varrho}$}
\end{equation}
which implies that 
\begin{equation}
  |\w{R^wu,u} | \le 2{\varrho}\left(\mn{F_h^wu}^2 +
    h^{-\frac{2}{k+1}}\mn{u}^2\right) 
\end{equation}
for $u \in C_0^\infty$ and  small $h >0$. By~\eqref{pouest}
and~\eqref{errest} we find
\begin{gather*}
\sum_{j} \mn{{\phi}_j^w(x,D_x)(D_t + i
  F_{h,w_j}(t))u}^2 \le \\ 4 \mn{(D_t + i
  F_{h}(t))u}^2 +  8{\varrho}\left(\mn{F_h^wu}^2 +
    h^{-\frac{2}{k+1}}\mn{u}^2\right)
\end{gather*}
for $u \in C_0^\infty$  and  small enough $h >0$.
We obtain from ~\eqref{pouest} and~\eqref{locest} that 
\begin{gather*}
 h^{-\frac{2}{k+1}}\mn{u}^2 + \mn{F^w_hu}^2  
\le \\ 16C^2\mn{(D_t +
     iF_h^w)u}^2  +  8{\varrho}(4C^2 + 1)\left(\mn{F_h^wu}^2 +
     h^{-\frac{2}{k+1}}\mn{u}^2\right)
\end{gather*}
for $u \in C_0^\infty$   and  small $h >0$.
Thus for small enough ${\varrho}$ and $h$ we obtain the estimate
\begin{equation*}
 h^{-\frac{2}{k+1}}\mn{u}^2 \le 32C^2\mn{(D_t +
     iF_h^w)u}^2 = 32C^2\mn{P_hu}^2 \qquad u \in C_0^\infty.
\end{equation*}
Since $P(h) = hP_h$ we obtain~\eqref{propest}, which completes the
proof of the Proposition.
\end{proof}

\section{Dissipative operators without the dynamical condition}
\label{bn}

The failure of \eqref{eq:dyn} and the consequent failure of Theorem \ref{thm2}
are illustrated by the following variant of Davies's  example \cite{Dav1}
in dimension two:
\[  p ( x , \xi ) = \xi_1^2 + \xi_2^2 + x_1^2 - i x_2^2  \,.\]
The spectrum of $ p^w ( x , h D) $  is given by $ \{ ( 2 n + 1 ) h \} _{
n \in \NN} \cup \{ e^{ \pi i / 4 } ( 2k + 1 )h  \}_{ k \in \NN } $, and 
$ \Lambda ( p ) $ is the first quadrant. We see that \eqref{eq:cone} is 
satisfied everywhere at the boundary, \eqref{eq:prin} everywhere except
for $ z_0= 0 $. The dynamical condition \eqref{eq:dyn} is satisfied
on the imaginary half-axis but fails on the real half-axis.

In this section we present a general result in the same spirit.
Suppose that $ P ( h ) $ satisfies the general assumptions of \S \ref{rsq},
but does not need to have an analytic symbol.
In 
addition we assume that
\begin{gather}
\label{eq:strong}
\begin{gathered}
P( h ) = Q ( h ) - i  W ( h ) \\
Q ( h ) = Q ( h )^* \,, \ \ W ( h ) \geq 0 \,. 
\end{gathered}
\end{gather}
In the classical terminology of 
\cite{gk} this means that our non-self-adjoint operator is dissipative.
One way of  achieving this semi-classically is by putting 
\begin{gather*}
Q ( h ) = q^w ( x , h D) \,, \ \text{ $ q $ real valued,} \ \ \ 
W ( h ) = a^{\rm{W}} ( x , h D) \,, \ \text{ with $ a \geq 0 $}\,, \\
a^{\rm{W}} ( x , h D) = \int \! \! \! 
\int a ( y , \eta) \Gamma^w_{ y , \eta} ( x ,
h D ) dy d \eta \,, \ \
\Gamma_{ y , \eta} ( x, \xi ) = \frac{1}{ \pi^n } e^{ - |x - y |^2 - 
| \xi - \eta |^2 } \,. 
\end{gather*}
Here $ \bullet^w $ stands for the Weyl quantization \eqref{eq:2.weyl},
and $ \bullet^{\rm{W}} $ for the Wick quantization. 

This assumption immediately implies that $ \Lambda ( p ) \subset \{\Im z \leq 0 
\} $, and also that
\[ z \in \sigma ( P ( h ))  \ \Longrightarrow 
\Im z \leq 0 \,, \]
In fact we have
\begin{gather}
\label{eq:3.1}
\begin{gathered}
- \Im \langle ( P ( h ) - z ) u , u \rangle \geq \Im z \| u \|^2 \\
( P ( h ) - z )^{-1} = {\mathcal O} ( 1/ \Im z ) \,, \ \ \Im z > 0 \,.
\end{gathered}
\end{gather}

We can now use techniques common in the study of dissipative operators
-- see \cite{gk},\cite{Mar},\cite{Hit}, and references given there. Similar 
techniques have been also developed in the study of semi-classical
resonances \cite{Sjn},\cite{StVo},\cite{St},\cite{TaZw} and our approach
follows these works in an easier setting of dissipative operators.

Thus we start with 
\begin{lem}
\label{p:3.1}
Assuming \eqref{eq:strong} 
we have, 
 for any open and precompact subset, $ \Omega $,  of any
component of $ \CC \setminus \Sigma_\infty (p) $ intersecting $ \complement 
\Sigma ( p ) $,
\begin{equation}
\label{eq:p.3}
\| ( P ( h ) - z )^{-1} \| \leq \exp \left(C_\Omega
 h^{-n} \log \frac1{g ( h ) }
\right)  \,, \ \ 
z \in \Omega \setminus \bigcup_{ z_j \in \sigma ( P ( h )) } D( z_j , g( h ) )
\,.
\end{equation}
\end{lem}
\begin{proof}
We can assume, without loss of generality, that $ \Omega = D( z_0 , 
\epsilon ) $, and that $ D( z_0 , 3 \epsilon) \cap \Sigma_\infty ( p ) = 
\emptyset $. It then follows that for $ C $ sufficiently large
$  p ( x , \xi ) \notin D ( z_0 , 2 \epsilon ) 
$, and $ | ( x , \xi ) | > C $. 

We can then
find $ p^\# \in \CI_{\rm{b}} ( T^* \RR^n ) $ such that 
\begin{gather*} p^\# ( x , \xi ) = p ( x , \xi )  \ \ \text{ for $ | ( x , 
\xi ) | > C $,} \\
\forall \;  z \in \Omega  \ ( x, \xi ) \in T^* \RR^n \ \ 
|  ( p ^\#  ( x , \xi ) - z ) ^{ -1} | \leq C \,. 
\end{gather*}
If fact, choose 
 $ \alpha \; : \; \CC  \rightarrow \CC \setminus D ( 0 , 2 \epsilon)  $, such
that $ \alpha ( w ) = w $ 
on $ \CC \setminus D ( z_0 , 3 \epsilon ) $, and put $ p^\# = \alpha \circ
p $. 

This shows that for $ h $ small enough $ p^\# ( x, h D ) - z $ is 
invertible, $ z \in \Omega $. 

In view of the compact support of the symbol $ p - p^\# $ we have 
\begin{gather*}
  p ( x , h D) - p^{\#} ( x , h D ) = A + B \\
A \; : \; L^2 ( \RR^n ) \ 
\longrightarrow \ \CI_{\rm{c}} ( \RR^n ) \,, \\
B = {\mathcal O} ( h^\infty ) \; : \; L^2 ( \RR^n ) \ \longrightarrow \
L^2 ( \RR^n ) \  \,, \\
p ( x , h D) - z = ( p^\# ( x , h D ) + B - z ) ( I + ( P^\# + B - z ) ^{-1}
A ) \,.
\end{gather*}
Hence  $ (  p^\# ( x, h D ) - z ) ^{-1} ( p ( x , h D) - 
p^{\#} ( x , h D ) )$ is a sum of a compact operator on $ L^2 ( \RR^n ) $
and an operator of small norm. 

Since the operator $ K ( z ) = 
(  p^\# ( x, h D ) + B - z ) ^{-1}  A   $ is
compact we can define
\[ g ( z) = \det ( I + K ( z ) ) \,,\]
which is holomorphic in $ \Omega $. The inverse at a point $ z_0
\notin \Sigma ( p ) $, 
$ ( I + K ( z_0) )^{-1} $, exists for $ h $ sufficiently small and it 
is bounded independently of $ h $.  We also see that
\[ |\det ( I + K ( z_0 )) |^{-1}  \leq C_1 \,, \ \text{with $ C_1 $ 
independent of $ h$.} \]
As in Proposition \ref{p:3}, $ ( I + K ( z) )^{-1} $ is meromorphic
and, following \cite[Ch.5, Theorem 3.1]{gk}, we see that
\[ \|  ( I + K ( z) )^{-1} \| \leq \frac{ \det
( I + ( K ( z) K(z)^*)^{\frac12} )}{ |\det ( I + K ( z ) )| } \,.\]
Hence we need to estimate the determinants from above and below.
The upper bound is clear from $ |\det ( 1 + A ) | \leq \exp \tr (AA^*)^{
\frac12} $: it is given by $ \exp  {\mathcal O} ( h^{-n} )$ since
\[ \tr \psi ( x ) b^w ( x , h D ) \langle h D\rangle
^{-m}  = \frac{1}{( 2 \pi h )^n}  \int
\int \psi ( x) a ( x , \xi ) \langle \xi \rangle^{-m} dx d\xi 
= {\mathcal  O} ( h^{-n} ) \,. \]
The zeros of $ \det ( I + K ( z ) ) $ coincide with the eigenvalues
of $ P ( h ) $, and
the usual complex analytic methods (see for instance \cite[Ch.1]{Mar},
and \cite{Sjn}) show that 
\[ |\det ( I + K ( z) ) |^{-1} \leq \frac{1}{g (h)}  \exp ( C h^{-n} ) 
\,,  \ \ z \in \Omega  \setminus 
\bigcup_{ z_j \in \sigma ( P ( h )) } D( z_j , g ( h ) ) \,, \]
which proves the lemma.
\end{proof}

The main result of this section will be an application of 
Lemma \ref{p:3.1}, \eqref{eq:3.1} and of the following simple
function theoretical lemma similar to \cite[Lemma 2]{TaZw}:
\begin{lem}
\label{p:3.2}
Suppose that $ F ( z ) $ is holomorphic in $ [ - \delta , \delta ] 
+ i [ - \epsilon , \epsilon ] $, and $ | F ( z ) | \leq M$, $ M \geq 2 $,
there.
Suppose in addition that
\[ | F ( z ) | \leq \frac{1}{ \Im z } \ \ \text{ for $ \Im z > 0 $}\,, \ \ 
\text{ and that $ \ \ \frac{\epsilon}{\delta} \ll \log M $.} \]
Then 
\[ | F ( z ) | \leq 2 \frac{\log M } {\epsilon} \,, \ \ \text{ for 
$ | z | \leq \delta/2 $ and $ \Im z = 0 $.} \]
\end{lem}
\begin{proof} 
The assumption on $ \delta/ \epsilon $ allows us to construct a
holomorphic function, $ u ( z ) $,  such that $ |u ( z ) | > 2/3 $
for $ \Im z = 0$, and $ |z|  < \delta/2 $, and $ | u ( z ) | \ll 
1/ M $ for $ |\Re z | = \delta $ and $ |\Im z | \leq \epsilon $.
If we apply an optimized
``three line theorem'' argument to $ u ( z ) F ( z ) $
the lemma follows. 
\end{proof}

The two lemmas immediately give
\begin{prop}
\label{t:3}
Suppose that $ P ( h ) $ satisfies \eqref{eq:strong}. For $ E \in \RR$ and 
$ k > 0 $ put
\begin{gather*} \Omega
( h ) = [ E -  \delta(h) , E +  \delta(h)  ] - i [0, - K 
h^{ n} \delta(h) ] \,, \ \ \text{$K$ large and fixed.}
\end{gather*}
Then $ h $ small enough we have
\begin{equation}
\label{eq:p.5} 
 \sigma ( P ( h ) ) \cap \Omega ( h ) = \emptyset \
\Longrightarrow \ 
\| ( P ( h ) - \lambda )^{-1} \| \leq C h^{- 2 n  } \delta(h)^{-1} \,, 
\ z \in \widetilde \Omega ( h ) \,. \end{equation}
\end{prop}

This means that in small neighbourhoods of the real axis, only 
eigenvalues can produce extreme growth of the resolvent.
A contradiction argument from 
\cite{St},\cite{StVo},\cite{TaZw} now shows that existence of quasi-modes
for the operator $ P ( h )$ implies the existence of spectrum arbitrarily 
close to the real axis, which depending on the structure of $ W ( h )$,
can be the boundary of $ \Sigma ( P ) $.  

In particular we have
\begin{prop}
\label{t:4}
Suppose that $ P ( h ) $ satisfies \eqref{eq:strong} and that 
there exists $ u ( h ) \neq 0 $ and $ \lambda ( h )\in \RR $ such that 
\begin{equation}
\label{eq:qm} \| ( P ( h )  - \lambda ( h ) ) u (h ) \| = {\mathcal O} ( h^N ) \| u ( h ) \|  \,,
\end{equation}
for some large $ N$. Then for $ h $ small
enough, and with $ K $ sufficiently large 
\[  d ( \sigma ( P ( h ) ) , \lambda ( h ) ) < K h^{N -n }  \,. 
\]
\end{prop}
This is clear from \eqref{eq:p.5} since applying that resolvent estimate to \eqref{eq:qm}
we obtain a contradiction.

In the analytic case we expect that if $ \Im p $ vanishes to a high order on a closed
orbit of $ H_{\Re p} $ then we can construct a quasi-mode $ u ( h ) $ satisfying \eqref{eq:qm}
with a fixed $ N$. If we allow $ \CI $ coefficients then constructions of quasi-modes 
with arbitrarily large $ N $s are possible for operators satisfying the assumptions of this 
section -- see \cite{St},\cite{TaZw} and references given there, and also \cite[Fig.1]{Zw1}
for a figure of an example.

\end{document}